\documentclass[11pt]{amsart}
\usepackage{amsmath, amssymb, latexsym}
\usepackage{mathrsfs}
\usepackage[all, knot]{xy}
\usepackage{bm}
\usepackage{color, url}
\xyoption{arc}

\theoremstyle{definition}

\numberwithin{equation}{section}

\numberwithin{equation}{section}

\newtheorem{theorem}{Theorem}[section]
\newtheorem{corollary}[theorem]{Corollary}
\newtheorem{lemma}[theorem]{Lemma}
\newtheorem{proposition}[theorem]{Proposition}
\newtheorem{conjecture}[theorem]{Conjecture}
\newtheorem{example}[theorem]{Example}

\theoremstyle{definition}
\newtheorem{definition}[theorem]{Definition}
\newtheorem{remark}[theorem]{Remark}

\newcommand{\HOM}{\text{HOM}}
\newcommand{\Hom}{\text{Hom}}

\newcommand{\B}{\mathbb{B}}

\newcommand{\E}{\mathcal{E}}
\newcommand{\F}{\mathcal{F}}
\newcommand{\Q}{\mathbb{Q}}
\newcommand{\R}{\mathcal{R}}
\newcommand{\Z}{\mathbb{Z}}
\newcommand{\N}{\mathbb{N}}
\newcommand{\h}{\mathfrak{h}}
\newcommand{\g}{\mathfrak{g}}
\newcommand{\K}{\mathbb{K}}

\newcommand{\ii}{\textit{\textbf{i}}}
\newcommand{\jj}{\textit{\textbf{j}}}
\newcommand{\kk}{\textit{\textbf{k}}}
\newcommand{\iii}{\textbf{i}}
\newcommand{\aaa}{\textbf{a}}

\newcommand{\Pol}{\mathscr{P}}

\newcommand{\Ind}{\text{Ind}}
\newcommand{\Res}{\text{Res}}

\newcommand{\gdim}{\mathbf{Dim}}
\newcommand{\Mod}{\text{Mod}}
\newcommand{\pMod}{\text{pMod}}
\newcommand{\fMod}{\text{fMod}}

\newcommand{\Ch}{\text{Ch}}

\newcommand{\hd}{\text{hd}}

\newcommand{\Max}{\text{max}}

\newcommand{\Spr}{\text{Spr}}

\newcommand{\genO}[3] 
{
\fontsize{9}{9}\selectfont
\xy
(0,5)*{}; (0,-5)*{} **\dir{-};
(5,0)*{\cdots};
(10,5)*{}; (10,-5)*{} **\dir{-};
(15,0)*{\cdots};
(20,5)*{}; (20,-5)*{} **\dir{-};
(0,-7)*{#1}; (10,-7)*{#2}; (20,-7)*{#3};
\endxy
\fontsize{10}{10}\selectfont}
\newcommand{\genX}[3] 
{
\fontsize{9}{9}\selectfont
\xy
(0,5)*{}; (0,-5)*{} **\dir{-};
(5,0)*{\cdots};
(10,5)*{}; (10,-5)*{} **\dir{-};
(10,0)*{\bullet}; (15,0)*{\cdots};
(20,5)*{}; (20,-5)*{} **\dir{-};
(0,-7)*{#1}; (10,-7)*{#2}; (20,-7)*{#3};
\endxy
\fontsize{10}{10}\selectfont}
\newcommand{\genT}[4] 
{\fontsize{8}{8}\selectfont
\xy
(0,5)*{}; (0,-5)*{} **\dir{-};
(5,0)*{\cdots};
(11,5)*{}; (20,-5)*{} **\dir{-};
(20,5)*{}; (11,-5)*{} **\dir{-};
(27,0)*{\cdots};
(34,5)*{}; (34,-5)*{} **\dir{-};
(0,-7)*{#1}; (9.5,-7)*{#2}; (23,-7)*{#3}; (36,-7)*{#4};
\endxy
\fontsize{10}{10}\selectfont}

\newcommand{\dcross}[2]
{\fontsize{9}{9}\selectfont
\xy
(0,7)*{}="T1"; (7,7)*{}="T2";
(0,-7)*{}="B1"; (7,-7)*{}="B2";
"T1"; "B1" **\crv{(11, 0)};
"T2"; "B2" **\crv{(-4,0)};
(-1,-9)*{#1}; (8,-9)*{#2};
\endxy
\fontsize{10}{10}\selectfont}
\newcommand{\dcrossA}[2]
{\fontsize{9}{9}\selectfont
\xy
(0,5)*{}; (0,-5)*{} **\dir{-};
(6,5)*{}; (6,-5)*{} **\dir{-};
(-1,-7)*{#1}; (7,-7)*{#2};
\endxy
\fontsize{10}{10}\selectfont}
\newcommand{\dcrossAA}[2]
{\fontsize{9}{9}\selectfont
\xy
(0,5)*{}; (0,-5)*{} **\dir{-};
(6,5)*{}; (6,-5)*{} **\dir{-};
(0,-7)*{#1}; (6,-7)*{#2};
\endxy
\fontsize{10}{10}\selectfont}
\newcommand{\dcrossL}[3]
{ \fontsize{10}{10}\selectfont
\xy
(0,5)*{}; (0,-5)*{} **\dir{-};
(6,5)*{}; (6,-5)*{} **\dir{-};
(0,0)*{\bullet}; (-6,0)*{#1}; (-1,-7)*{#2}; (7,-7)*{#3};
\endxy
\fontsize{10}{10}\selectfont}
\newcommand{\dcrossR}[3]
{\fontsize{10}{10}\selectfont
\xy
(0,5)*{}; (0,-5)*{} **\dir{-};
(6,5)*{}; (6,-5)*{} **\dir{-};
(6,0)*{\bullet}; (12,0)*{#1}; (-1,-7)*{#2}; (7,-7)*{#3};
\endxy
\fontsize{10}{10}\selectfont}

\newcommand{\dcrossLL}[3]
{ \fontsize{9}{9}\selectfont
\xy
(0,5)*{}; (0,-5)*{} **\dir{-};
(6,5)*{}; (6,-5)*{} **\dir{-};
(0,0)*{\bullet}; (-6,0)*{#1}; (-1,-7)*{#2}; (7,-7)*{#3};
\endxy
\fontsize{10}{10}\selectfont}
\newcommand{\dcrossRR}[3]
{\fontsize{9}{9}\selectfont
\xy
(0,5)*{}; (0,-5)*{} **\dir{-};
(6,5)*{}; (6,-5)*{} **\dir{-};
(6,0)*{\bullet}; (12,0)*{#1}; (-1,-7)*{#2}; (7,-7)*{#3};
\endxy
\fontsize{10}{10}\selectfont}

\newcommand{\LU}[3]
{\fontsize{9}{9}\selectfont
\xy
(0,5)*{}; (8,-5)*{} **\dir{-} \POS?(.25)="x";
(8,5)*{}; (0,-5)*{} **\dir{-};
"x"*{\bullet}; "x"+(-2,2)*{#1}; (0,-7)*{#2}; (8,-7)*{#3}; \endxy
\fontsize{10}{10}\selectfont}
\newcommand{\RD}[3]
{\fontsize{9}{9}\selectfont
\xy
(0,5)*{}; (8,-5)*{} **\dir{-} \POS?(.75)="x";
(8,5)*{}; (0,-5)*{} **\dir{-};
"x"*{\bullet}; "x"+(2,2)*{#1}; (0,-7)*{#2}; (8,-7)*{#3}; \endxy
\fontsize{10}{10}\selectfont}
\newcommand{\RU}[3]
{\fontsize{9}{9}\selectfont
\xy
(0,5)*{}; (8,-5)*{} **\dir{-};
(8,5)*{}; (0,-5)*{} **\dir{-} \POS?(.25)="x";
"x"*{\bullet}; "x"+(2,2)*{#1}; (0,-7)*{#2}; (8,-7)*{#3}; \endxy
\fontsize{10}{10}\selectfont}
\newcommand{\LD}[3]
{\fontsize{9}{9}\selectfont
\xy
(0,5)*{}; (8,-5)*{} **\dir{-};
(8,5)*{}; (0,-5)*{} **\dir{-} \POS?(.75)="x";
"x"*{\bullet}; "x"+(-2,2)*{#1}; (0,-7)*{#2}; (8,-7)*{#3}; \endxy
\fontsize{10}{10}\selectfont}

\newcommand{\BraidLLL}[3]
{\fontsize{9}{9}\selectfont
\xy
(0,7)*{}="T1"; (6,7)*{}="T2";  (12,7)*{}="T3";
(0,-7)*{}="B1"; (6,-7)*{}="B2"; (12,-7)*{}="B3";
"T1"; "B3" **\crv{(0,0)&(12,0)}; "T3"; "B1" **\crv{(12,0)&(0,0)};
"T2"; "B2" **\crv{(6,6)&(3,4.5)&(0,1.5)&(0,-1.5)&(3,-4.5)&(6,-6)};
(0,-9)*{#1}; (6,-9)*{#2};(12,-9)*{#3};
\endxy
\fontsize{10}{10}\selectfont}
\newcommand{\BraidRRR}[3]
{\fontsize{9}{9}\selectfont
\xy
(0,7)*{}="T1"; (6,7)*{}="T2";  (12,7)*{}="T3";
(0,-7)*{}="B1"; (6,-7)*{}="B2"; (12,-7)*{}="B3";
"T1"; "B3" **\crv{(0,0)&(12,0)}; "T3"; "B1" **\crv{(12,0)&(0,0)};
"T2"; "B2" **\crv{(6,6)&(8,4.5)&(12,1.5)&(12,-1.5)&(8,-4.5)&(6,-6)};
(0,-9)*{#1}; (6,-9)*{#2};(12,-9)*{#3};
\endxy
\fontsize{10}{10}\selectfont}

\newcommand{\BraidL}[3]
{\fontsize{8}{8}\selectfont
\xy
(0,7)*{}="T1"; (6,7)*{}="T2";  (12,7)*{}="T3";
(0,-7)*{}="B1"; (6,-7)*{}="B2"; (12,-7)*{}="B3";
"T1"; "B3" **\crv{(0,0)&(12,0)}; "T3"; "B1" **\crv{(12,0)&(0,0)};
"T2"; "B2" **\crv{(6,6)&(3,4.5)&(0,1.5)&(0,-1.5)&(3,-4.5)&(6,-6)};
(-1,-9)*{#1}; (6,-9)*{#2};(14,-9)*{#3};
\endxy
\fontsize{10}{10}\selectfont}
\newcommand{\BraidR}[3]
{\fontsize{8}{8}\selectfont
\xy
(0,7)*{}="T1"; (6,7)*{}="T2";  (12,7)*{}="T3";
(0,-7)*{}="B1"; (6,-7)*{}="B2"; (12,-7)*{}="B3";
"T1"; "B3" **\crv{(0,0)&(12,0)}; "T3"; "B1" **\crv{(12,0)&(0,0)};
"T2"; "B2" **\crv{(6,6)&(8,4.5)&(12,1.5)&(12,-1.5)&(8,-4.5)&(6,-6)};
(-1,-9)*{#1}; (6,-9)*{#2};(14,-9)*{#3};
\endxy
\fontsize{10}{10}\selectfont}
\newcommand{\threeDotStrands}[3]
{ \fontsize{9}{9}\selectfont
\xy
(0,7)*{}; (0,-7)*{} **\dir{-};
(6,7)*{}; (6,-7)*{} **\dir{-};
(12,7)*{}; (12,-7)*{} **\dir{-};
(0,0)*{\bullet};(12,0)*{\bullet}; (-3,0)*{c};(23,0)*{-\ell a_{ij}-1-c}; (0,-9)*{#1}; (6,-9)*{#2};(12,-9)*{#3};
\endxy
\fontsize{10}{10}\selectfont}

\newcommand{\dcrossI}[2]
{\fontsize{9}{9}\selectfont
\xy
(0,6)*{}="T1"; (6,6)*{}="T2";
(0,-6)*{}="B1"; (6,-6)*{}="B2";
"T1"; "B1" **\crv{(10, 0)};
"T2"; "B2" **\crv{(-4,0)};
(0,-8)*{#1}; (6,-8)*{#2};
\endxy
\fontsize{10}{10}\selectfont}
\newcommand{\BraidLI}[3]
{\fontsize{9}{9}\selectfont
\xy
(0,6)*{}="T1"; (5,6)*{}="T2";  (10,6)*{}="T3";
(0,-6)*{}="B1"; (5,-6)*{}="B2"; (10,-6)*{}="B3";
"T1"; "B3" **\crv{(0,0)&(10,0)}; "T3"; "B1" **\crv{(10,0)&(0,0)};
"T2"; "B2" **\crv{(5,5)&(3,4)&(0,1.5)&(0,-1.5)&(3,-4)&(5,-5)};
(0,-8)*{#1}; (5,-8)*{#2};(10,-8)*{#3};
\endxy
\fontsize{10}{10}\selectfont}
\newcommand{\BraidRI}[3]
{\fontsize{9}{9}\selectfont
\xy
(0,6)*{}="T1"; (5,6)*{}="T2";  (10,6)*{}="T3";
(0,-6)*{}="B1"; (5,-6)*{}="B2"; (10,-6)*{}="B3";
"T1"; "B3" **\crv{(0,0)&(10,0)}; "T3"; "B1" **\crv{(10,0)&(0,0)};
"T2"; "B2" **\crv{(5,5)&(7,4)&(10,1.5)&(10,-1.5)&(7,-4)&(5,-5)};
(0,-8)*{#1}; (5,-8)*{#2};(10,-8)*{#3};
\endxy
\fontsize{10}{10}\selectfont}

\newcommand{\BraidLII}[3]
{\fontsize{8}{8}\selectfont
\xy
(0,6)*{}="T1"; (5,6)*{}="T2";  (10,6)*{}="T3";
(0,-6)*{}="B1"; (5,-6)*{}="B2"; (10,-6)*{}="B3";
"T1"; "B3" **\crv{(0,0)&(10,0)}; "T3"; "B1" **\crv{(10,0)&(0,0)};
"T2"; "B2" **\crv{(5,5)&(3,4)&(0,1.5)&(0,-1.5)&(3,-4)&(5,-5)};
(-1,-8)*{#1}; (5,-8)*{#2};(11,-8)*{#3};
\endxy
\fontsize{10}{10}\selectfont}
\newcommand{\BraidRII}[3]
{\fontsize{8}{8}\selectfont
\xy
(0,6)*{}="T1"; (5,6)*{}="T2";  (10,6)*{}="T3";
(0,-6)*{}="B1"; (5,-6)*{}="B2"; (10,-6)*{}="B3";
"T1"; "B3" **\crv{(0,0)&(10,0)}; "T3"; "B1" **\crv{(10,0)&(0,0)};
"T2"; "B2" **\crv{(5,5)&(7,4)&(10,1.5)&(10,-1.5)&(7,-4)&(5,-5)};
(-1,-8)*{#1}; (5,-8)*{#2};(11,-8)*{#3};
\endxy
\fontsize{10}{10}\selectfont}
\newcommand{\dcrossII}[2]
{\fontsize{9}{9}\selectfont
\xy
(0,6)*{}="T1"; (6,6)*{}="T2";
(0,-6)*{}="B1"; (6,-6)*{}="B2";
"T1"; "B1" **\crv{(10, 0)};
"T2"; "B2" **\crv{(-4,0)};
(-1,-8)*{#1}; (7,-8)*{#2};
\endxy
\fontsize{10}{10}\selectfont}

\newcommand{\Down}
{\fontsize{10}{10}\selectfont
\xy
(0,0)*{}; (12,12)*{} **\dir{-};
(3,0)*{}; (15,12)*{} **\dir{-};
(6,0)*{}; (18,12)*{} **\dir{-};
(9,0)*{}; (0,12)*{} **\dir{-};
(12,0)*{}; (3,12)*{} **\dir{-};
(15,0)*{}; (6,12)*{} **\dir{-};
(18,0)*{}; (9,12)*{} **\dir{-};
(0,12)*{}; (9,12)*{} **\dir{-};
(0,17)*{}; (9,17)*{} **\dir{-};
(0,12)*{}; (0,17)*{} **\dir{-};
(9,12)*{}; (9,17)*{} **\dir{-};
(12,12)*{}; (18,12)*{} **\dir{-};
(12,17)*{}; (18,17)*{} **\dir{-};
(12,12)*{}; (12,17)*{} **\dir{-};
(18,12)*{}; (18,17)*{} **\dir{-};
(4.5, 14)*{e_{i,\ell}};(0,-2)*{j}; (3,-2)*{j};(6,-2)*{j};
(15, 14)*{e_{j,k}};(9,-2)*{i}; (12,-2)*{i};(15,-2)*{i};(18,-2)*{i};(25,7)*{=};
\endxy
\fontsize{11}{10}\selectfont}

\newcommand{\Up}
{\fontsize{10}{10}\selectfont
\xy
(0,2)*{}; (12,14)*{} **\dir{-};
(3,2)*{}; (15,14)*{} **\dir{-};
(6,2)*{}; (18,14)*{} **\dir{-};
(9,2)*{}; (0,14)*{} **\dir{-};
(12,2)*{}; (3,14)*{} **\dir{-};
(15,2)*{}; (6,14)*{} **\dir{-};
(18,2)*{}; (9,14)*{} **\dir{-};
(0,-3)*{}; (6,-3)*{} **\dir{-};
(0,2)*{}; (6,2)*{} **\dir{-};
(0,-3)*{}; (0,2)*{} **\dir{-};
(6,-3)*{}; (6,2)*{} **\dir{-};
(9,-3)*{}; (18,-3)*{} **\dir{-};
(9,2)*{}; (18,2)*{} **\dir{-};
(9,-3)*{}; (9,2)*{} **\dir{-};
(18,-3)*{}; (18,2)*{} **\dir{-};
(3, -1)*{e_{j,k}};(0,16)*{i}; (3,16)*{i};(6,16)*{i};
(13.5, -1)*{e_{i,\ell}};(9,16)*{i}; (12,16)*{j};(15,16)*{j};(18,16)*{j};
\endxy
\fontsize{11}{10}\selectfont}
\newcommand{\DDown}
{\fontsize{10}{10}\selectfont
\xy
(0,0)*{}; (12,12)*{} **\dir{-};
(3,0)*{}; (15,12)*{} **\dir{-};
(6,0)*{}; (18,12)*{} **\dir{-};
(9,0)*{}; (0,12)*{} **\dir{-};
(12,0)*{}; (3,12)*{} **\dir{-};
(15,0)*{}; (6,12)*{} **\dir{-};
(18,0)*{}; (9,12)*{} **\dir{-};
(0,-2)*{j}; (3,-2)*{j};(6,-2)*{j};
(9,-2)*{i}; (12,-2)*{i};(15,-2)*{i};(18,-2)*{i};
\endxy
\fontsize{11}{10}\selectfont}
\newcommand{\coset}
{\fontsize{10}{10}\selectfont
\xy
(-9,0)*{}; (-9,12)*{} **\dir{-};
(-6,0)*{}; (-6,12)*{} **\dir{-};
(-3,0)*{}; (-3,12)*{} **\dir{-};
(0,0)*{}; (9,12)*{} **\dir{-};
(3,0)*{}; (12,12)*{} **\dir{-};
(6,0)*{}; (15,12)*{} **\dir{-};
(9,0)*{}; (0,12)*{} **\dir{-};
(12,0)*{}; (3,12)*{} **\dir{-};
(15,0)*{}; (6,12)*{} **\dir{-};
(18,0)*{}; (18,12)*{} **\dir{-};
(21,0)*{}; (21,12)*{} **\dir{-};
(-1.5,-3)*{\underbrace{ \quad\quad\quad\quad \ }_{n}};
(15,-3)*{\underbrace{ \quad\quad\quad \  }_{\ell}};
(12,15)*{\overbrace{ \quad \  }^{k}};
\endxy
\fontsize{11}{10}\selectfont}

\newcommand{\genTT}[4] 
{\fontsize{9}{9}\selectfont
\xy
(0,5)*{}; (0,-5)*{} **\dir{-};
(4,0)*{\cdots};
(7,5)*{}; (13,-5)*{} **\dir{-};
(13,5)*{}; (7,-5)*{} **\dir{-};
(14,0)*{\cdots};
(20,5)*{}; (20,-5)*{} **\dir{-};
(0,-7)*{#1}; (7,-7)*{#2}; (14,-7)*{#3}; (20,-7)*{#4};
\endxy
\fontsize{10}{10}\selectfont}

\newcommand{\dcrossS}[2]
{\fontsize{9}{9}\selectfont
\xy
(0,7)*{}="T1"; (7,7)*{}="T2";
(0,-7)*{}="B1"; (7,-7)*{}="B2";
"T1"; "B1" **\crv{(11, 0)};
"T2"; "B2" **\crv{(-4,0)};
(0,-9)*{#1}; (7,-9)*{#2};
\endxy
\fontsize{10}{10}\selectfont}
\newcommand{\dcrossAS}[2]
{\fontsize{9}{9}\selectfont
\xy
(0,5)*{}; (0,-5)*{} **\dir{-};
(6,5)*{}; (6,-5)*{} **\dir{-};
(0,-7)*{#1}; (6,-7)*{#2};
\endxy
\fontsize{10}{10}\selectfont}
\newcommand{\dcrossLS}[3]
{ \fontsize{10}{10}\selectfont
\xy
(0,5)*{}; (0,-5)*{} **\dir{-};
(6,5)*{}; (6,-5)*{} **\dir{-};
(0,0)*{\bullet}; (-5,0)*{#1}; (0,-7)*{#2}; (6,-7)*{#3};
\endxy
\fontsize{10}{10}\selectfont}
\newcommand{\dcrossRS}[3]
{\fontsize{10}{10}\selectfont
\xy
(0,5)*{}; (0,-5)*{} **\dir{-};
(6,5)*{}; (6,-5)*{} **\dir{-};
(6,0)*{\bullet}; (11,0)*{#1}; (0,-7)*{#2}; (6,-7)*{#3};
\endxy
\fontsize{10}{10}\selectfont}

\newcommand{\LUS}[3]
{\fontsize{9}{9}\selectfont
\xy
(0,5)*{}; (8,-5)*{} **\dir{-} \POS?(.25)="x";
(8,5)*{}; (0,-5)*{} **\dir{-};
"x"*{\bullet}; "x"+(-2,2)*{#1}; (0,-7)*{#2}; (8,-7)*{#3}; \endxy
\fontsize{10}{10}\selectfont}
\newcommand{\RDS}[3]
{\fontsize{9}{9}\selectfont
\xy
(0,5)*{}; (8,-5)*{} **\dir{-} \POS?(.75)="x";
(8,5)*{}; (0,-5)*{} **\dir{-};
"x"*{\bullet}; "x"+(2,2)*{#1}; (0,-7)*{#2}; (8,-7)*{#3}; \endxy
\fontsize{10}{10}\selectfont}
\newcommand{\RUS}[3]
{\fontsize{9}{9}\selectfont
\xy
(0,5)*{}; (8,-5)*{} **\dir{-};
(8,5)*{}; (0,-5)*{} **\dir{-} \POS?(.25)="x";
"x"*{\bullet}; "x"+(2,2)*{#1}; (0,-7)*{#2}; (8,-7)*{#3}; \endxy
\fontsize{10}{10}\selectfont}
\newcommand{\LDS}[3]
{\fontsize{9}{9}\selectfont
\xy
(0,5)*{}; (8,-5)*{} **\dir{-};
(8,5)*{}; (0,-5)*{} **\dir{-} \POS?(.75)="x";
"x"*{\bullet}; "x"+(-2,2)*{#1}; (0,-7)*{#2}; (8,-7)*{#3}; \endxy
\fontsize{10}{10}\selectfont}

\newcommand{\BraidLS}[3]
{\fontsize{9}{9}\selectfont
\xy
(0,7)*{}="T1"; (6,7)*{}="T2";  (12,7)*{}="T3";
(0,-7)*{}="B1"; (6,-7)*{}="B2"; (12,-7)*{}="B3";
"T1"; "B3" **\crv{(0,0)&(12,0)}; "T3"; "B1" **\crv{(12,0)&(0,0)};
"T2"; "B2" **\crv{(6,6)&(3,4.5)&(0,1.5)&(0,-1.5)&(3,-4.5)&(6,-6)};
(0,-9)*{#1}; (6,-9)*{#2};(12,-9)*{#3};
\endxy
\fontsize{10}{10}\selectfont}
\newcommand{\BraidRS}[3]
{\fontsize{9}{9}\selectfont
\xy
(0,7)*{}="T1"; (6,7)*{}="T2";  (12,7)*{}="T3";
(0,-7)*{}="B1"; (6,-7)*{}="B2"; (12,-7)*{}="B3";
"T1"; "B3" **\crv{(0,0)&(12,0)}; "T3"; "B1" **\crv{(12,0)&(0,0)};
"T2"; "B2" **\crv{(6,6)&(8,4.5)&(12,1.5)&(12,-1.5)&(8,-4.5)&(6,-6)};
(0,-9)*{#1}; (6,-9)*{#2};(12,-9)*{#3};
\endxy
\fontsize{10}{10}\selectfont}
\newcommand{\threeDotStrandsS}[3]
{ \fontsize{9}{9}\selectfont
\xy
(0,7)*{}; (0,-7)*{} **\dir{-};
(5,7)*{}; (5,-7)*{} **\dir{-};
(10,7)*{}; (10,-7)*{} **\dir{-};
(0,0)*{\bullet};(10,0)*{\bullet}; (-3,0)*{c};(21,0)*{-a_{ij}-1-c}; (0,-9)*{#1}; (5,-9)*{#2};(10,-9)*{#3};
\endxy
\fontsize{10}{10}\selectfont}

\title[Categorification of quantum Borcherds-Bozec algebras]
{Categorification of quantum Borcherds-Bozec algebras}

\author[Stefano Kang]{Stefano V. Kang}
\address{Korea Research Institute of Arts and Mathematics,
Asan-si, Chungcheongnam-do, 31551, Korea}
\email{stefano.kang47@gmail.com}

\thanks{}

\author[Young Rock Kim]{Young Rock Kim${}^{*}$}
\address{Graduate School of Education, Hankuk University of Foreign Studies, Seoul, 02450,  Korea}
\email{rocky777@hufs.ac.kr} %
\thanks{${}^{*}$ Corresponding author.}
\author[Bolun Tong]{Bolun Tong}
\address{Graduate School of Education, Hankuk University of Foreign Studies, Seoul, 02450, Korea}
\email{tbl\_2018@hufs.ac.kr}



\begin{document}
\maketitle
\begin{abstract}
We categorify the quantum Borcherds-Bozec algebra $U_q(\g)$  for an arbitrary Borcherds-Cartan datum   by constructing their associated Khovanov-Lauda-Rouquier algebras.  In the Jordan quiver case, we show that the cyclotomic   Khovanov-Lauda-Rouquier algebras  provide a categorification of the irreducible highest weight $U_q(\g)$-modules.
\end{abstract}

\section*{\textbf{Introduction}}

\vskip 2mm

The \textit{Khovanov-Lauda-Rouquier algebras} (or \textit{quiver Hecke algebras}) were independently discovered by Khovanov-Lauda \cite{KL2009, KL2011} and Rouquier \cite{Rou}. In the Kac-Moody case, the category of finitely generated graded projective modules over Khovanov-Lauda-Rouquier algebras provides a categorification of the corresponding quantum groups, and for symmetric Cartan datum, the indecomposable projective modules correspond to Lusztig's canonical basis \cite{Rou1, VV2011}. The cyclotomic quotients of Khovanov-Lauda-Rouquier algebras categorify the irreducible highest weight representations of quantum groups and their crystals \cite{KK2012, LV2011}.

The \textit{quantum Borcherds-Bozec algebras} were introduced by T. Bozec \cite{Bozec2014b} in his study of perverse sheaves on quiver representation varieties, possibly with loops. He showed that the Grothendieck group arising from Lusztig sheaves is generated by the elementary simple perverse sheaves, answering a question posed by Lusztig in \cite{Lus93}.

A quantum Borcherds-Bozec algebra $ U_q(\mathfrak{g}) $ is determined by a Borcherds-Cartan datum, consisting of an index set $ I $ and a symmetrizable Borcherds-Cartan matrix $ A = (a_{ij})_{i,j \in I} $, where diagonal entries can be $ \leq 0 $. For an imaginary index $ i $ (i.e., $ a_{ii} \leq 0 $), there are infinitely many generators $ E_{i\ell}, F_{i\ell} $ ($ \ell \in \mathbb{Z}_{>0} $) associated with it.

In this paper, we apply Khovanov-Lauda's categorification theory to quantum Borcherds-Bozec algebras. To this end, for the imaginary indices,  we distinguish between the case where $a_{ii}=0$ (the Jordan quiver case) and the case where $a_{ii}<0$ (quiver with at least two loops). When $ a_{ii} = 0$, it is natural to relate the quantum Borcherds-Bozec algebra (which coincides with the classical Hall algebra) to the representations of symmetric groups. We construct the Khovanov-Lauda-Rouquier algebra as a deformation of the symmetric group algebras in this case. When $ a_{ii} < 0 $, higher-order comultiplications of $ F_{i\ell} $ present a problem. Thus, we use the \textit{primitive generators} $ \mathtt{b}_{i\ell} $, which have a simpler comultiplication, to provide our categorification. In this case, we treat each pair $ (i, \ell) $ as a simple root.

By considering the imaginary indices in this way and combining Khovanov-Lauda's work on Kac-Moody cases, we construct a Khovanov-Lauda-Rouquier type algebra $ R $ for quantum Borcherds-Bozec algebras using braid-like planar diagrams. Let $ K_0(R) $ be the Grothendieck group of the category of finitely generated graded projective $ R $-modules. We prove that there is a bialgebra isomorphism $ \Gamma$ between $ U^- $  and $K_0(R)$. Composing with an automorphism $\Psi $ of $ U^- $, we obtain $\mathcal{A}$-algebra isomorphism:

$$
\Phi = \Psi \Gamma^{-1}: K_0(R) \xrightarrow{\sim} {_\mathcal{A}} U^-,
$$
where $ \mathcal{A}= \mathbb{Z}[q, q^{-1}] $ and ${_\mathcal{A}} U^- $ is the $ \mathcal{A} $-form of $ U^- $. We conjecture the following:

\vskip 1mm
\noindent \textbf{Conjecture.}
Under the isomorphism $ \Phi $, the self-dual indecomposable projective modules of $ K_0(R) $ coincide with the elements of the canonical basis of $ {_\mathcal{A}} U^- $ given in \cite{Bozec2014b}.
\vskip 1mm

We verify this conjecture in the case of quivers with one vertex. For more general cases, we refer to \cite{VV2011}, which considers the Steinberg-type variety of a given quiver and provides a geometric realization of the Khovanov-Lauda-Rouquier algebras in Kac-Moody cases.

As an application of our construction of the Khovanov-Lauda-Rouquier algebra, we show in the Jordan quiver case that the cyclotomic algebra $ R^{\Lambda} $ $( \Lambda \in P^+ )$ provides a categorification of the irreducible highest weight module $ V(\Lambda)$. Essentially, we categorify the commutation relations of the generators $ E_{i\ell} $ and $ F_{it} $ when $ a_{ii} = 0 $.

We also consider a subalgebra of $ U^- $ which has a natural geometric interpretation as $ U^- $ (see Section 1.3). We construct the Khovanov-Lauda-Rouquier algebra for this subalgebra in the Appendix.

\vskip 3mm

\noindent\textbf{Acknowledgements.}
  Young Rock Kim and Bolun Tong were supported by the National Research Foundation of Korea (NRF) grant funded by the Korea government (MSIT) (No. 2021R1A2C1011467). Young Rock Kim was supported by Hankuk University of Foreign Studies Research Fund.


\vskip 6mm

\section{\textbf{Negative parts of quantum Borcherds-Bozec algebras}}
\vskip 3mm
\subsection{Notation}\
\vskip 3mm
In this paper, we fix an algebraically closed field $\K$ of characteristic zero.

Let $A$ be a $\Z$-graded $\K$-algebra.  For a graded $A$-module $M=\bigoplus_{n\in \Z}M_n$, its graded dimension is defined to be
$$\gdim \, M=\sum_{n\in\Z}(\text{dim}_{\K}M_n)q^n,$$
where $q$ is a formal variable. For $m\in \Z$,  the  degree shifted module
$M\{m\}$ is the graded $A$-module obtained from $M$ by putting $(M\{m\})_n=M_{n-m}$.  More generally, for  $f(q)=\sum_{m\in \Z}a_mq^m\in\N[q,q^{-1}]$, we set $M^{f}=\bigoplus_{m\in\Z}(M\{m\})^{\oplus a_m}$.

 Given two graded $A$-modules $M$ and $N$, we denote by $\text{Hom}_{A\text{-gr}}(M,N)$  the ${\K}$-vector space of  degree-preserving homomorphisms and form the $\Z$-graded vector space
$$\text{HOM}_A(M,N)=\bigoplus_{n\in\Z}\text{Hom}_{A\text{-gr}}(M\{n\},N)=\bigoplus_{n\in\Z}\text{Hom}_{A\text{-gr}}(M,N\{-n\}).$$
\vskip 3mm
We  use the term gr-projective (resp. gr-irreducible, gr-free and so on) module for the graded projective (resp. graded irreducible, graded free) module.
\vskip 3mm
For each $n\geq 0$, we denote $\lambda\vdash {n}$ ($\lambda\vDash{n}$) when $\lambda$  is a partition (composition) of $n$, and denote by $\mathcal P_n$ (resp. $\mathcal C_n$) the set of partitions (resp. compositions) of $n$.

\vskip 5mm
\subsection{Borcherds-Cartan datum and $U^-$}\

\vskip 3mm
Let $I$ be a finite index set. A Borcherds-Cartan datum $(I,A,\cdot)$ consists of
\begin{itemize}
\item [(a)] an integer-valued matrix $A=(a_{ij})_{i,j\in I}$ satisfying
\begin{itemize}
\item [(i)] $a_{ii}=2,0,-2,-4,\dots$,
\item [(ii)] $a_{ij}\in \Z_{\leq 0}$ for $i\neq j$,
\item [(iii)] there is a diagonal matrix $D=\text{diag}(r_i\in\Z_{>0}\mid i\in I)$ such that $DA$ is symmetric.
\end{itemize}
\item [(b)] a symmetric bilinear form $\alpha,\alpha'\mapsto \alpha\cdot \alpha'$ on $\Z[I]$ taking values in $\Z$, such that $$i\cdot j=r_ia_{ij}=r_ja_{ji} \ \ \text{for all} \ i,j\in I.$$
\end{itemize}

\vskip 3mm
 We set $I^{+}=\{ i \in I \mid a_{ii}=2 \}$, $I^{0}= \{ i \in I \mid a_{ii} =0 \}$,
$I^{-} =\{ i \in I \mid a_{ii}<0 \}$ and $I^{\leq 0} =I^0\cup I^-$. The elements in $I^+$ (resp. $I^{\leq 0}$) are called real indices (resp. imaginary indices).

\vskip 3mm
Let $q$ be an indeterminate. For each $i\in I$, we set
$$q_i=q^{r_i},\quad q_{(i)}=q^{\frac{i\cdot i}{2}}.$$
For each $i \in I^{+}$ and $n \in \N$, we set
$$[n]_i=\frac{q_i^n-q_i^{-n}}{q_i-q_i^{-1}}\ \text{and} \ [n]_i!=[n]_i[n-1]_i\cdots [1]_i.$$

\vskip 3mm
\begin{definition}\label{U-} Let $I^{\infty}= I^{+} \sqcup (I^{\leq 0}\times \Z_{>0})$. The negative part $U^- = U^-_q(\g)$ of the quantum Borcherds-Bozec algebra  associated with a given Borcherds-Cartan datum $(I, A, \cdot)$   is the associative algebra over $\Q(q)$  generated by $F_{i \ell}$ $((i,\ell)\in I^\infty)$,  satisfying the following relations
\begin{equation*}
\begin{aligned}
& \sum_{r+s=1-\ell a_{ij}}(-1)^r
{F_i}^{(r)}F_{j\ell}F_i^{(s)}=0 \quad \text{for} \ i\in
I^{+},(j,\ell)\in I^{\infty} \ \text {and} \ i \neq (j,\ell),\\
&F_{i\ell}F_{jk}-F_{jk}F_{i\ell} =0 \quad  \text{for}\  a_{ij}=0.
\end{aligned}
\end{equation*}
Here we denote $F_i^{(n)}=F_i^n /[n]_i!$ for $i\in I^+$ and $n \in \N$. The algebra $U^-$ is $\N[I]$-graded by assigning $|F_{i\ell}|=\ell i$.
\end{definition}

\vskip 3mm
Define a twisted multiplication on $U^-\otimes U^-$ by
$$(x_1\otimes x_2)(y_1\otimes y_2)=q^{-(|x_2|, |y_1|)}x_1y_1\otimes x_2y_2,$$
for homogeneous  $x_1,x_2,y_1,y_2$. By \cite[Proposition 14]{Bozec2014b}, we have an algebra homomorphism $\rho\colon U^-\rightarrow U^-\otimes U^-$ (with respect to the twisted multiplication on $U^-\otimes U^-$) given by $$\rho(F_{i\ell})=\sum_{m+n=\ell}q_{(i)}^{-mn}F_{im}\otimes F_{in} \ \ \text{for}  \ (i,\ell)\in I^{\infty},$$ and a nondegenerate symmetric bilinear form $\{ \ , \ \}\colon U^-\times U^-\rightarrow \Q(q)$ determined by
\begin{itemize}
\item[(i)] $\{x, y\} =0$ if $|x| \neq |y|$,
\item[(ii)] $\{1,1\} = 1$,
\item[(iii)] $\{F_{i\ell}, F_{i\ell}\}  \equiv 1  \pmod{q}$ for all $(i,\ell)\in I^\infty$,
\item[(iv)] $\{x, yz\} = \{\rho(x), y \otimes z\}$  for $x,y,z
\in U^-$.
\end{itemize}

\vskip 3mm
Let $\mathcal A=\Z[q,q^{-1}]$ be the ring of Laurent polynomials. The $\mathcal A$-form $_{\mathcal A} U^-$ is the $\mathcal A$-subalgebra  of $U^-$ generated by  $ F_i^{(n)}$ for $i\in I^+$, $ F_{i\ell}$ for $i\in I^{\leq 0},\ell\geq 1$.
\vskip 5mm
\subsection{Geometric setting for $U^-$ and  related algebras}\
\vskip 3mm

We briefly review the geometric construction for $U^-$ given in \cite{Bozec2014b}, \cite{Lus93} and \cite{LL09}.  Let $(I,H)$ be a quiver with vertices set $I$ and arrows set $H$. For each $h\in H$,   $h',h''\in I$ are   the origin
 and the goal of $h$ respectively. We allow    $h'$ and $h''$ to be equal.

Fix $\alpha=\sum_{i\in I}\alpha_i i\in \N[I]$. We set $V_\alpha=\bigoplus_{i\in I}\K^{\alpha_i}$, $E_\alpha=\bigoplus_{h\in H}\Hom(\K^{\alpha_{h'}},\K^{\alpha_{h''}})$ and $G_{\alpha}=\prod_{i\in I}\text{GL}_{\alpha_i}(\K)$. Let $G_{\alpha}$ acts on   $E_{\alpha}$ by $g\cdot(x_h)=(g_{h''}x_h g_{h'}^{-1})$.

 Denote by $D_{G_{\alpha}}(E_{\alpha})$  the  bounded $G_\alpha$-equivariant  derived category of $\K$-constructible
 complexes on $E_{\alpha}$ and by $P_{G_{\alpha}}(E_{\alpha})$ the abelian subcategory of $G_\alpha$-equivariant perverse sheaves.
\vskip 3mm

For a pair of sequences $\iii=(i_1,\dots,i_s)$ in $I$ and $\aaa=(a_1,\dots,a_s)$ in $\N$, we write $(\iii,\aaa)\vdash \alpha$ if $\alpha=\sum_{j}a_ji_j$. For such a pair, we set
$$ \mathscr F_{\iii,\aaa}=\{\text{all flags}\  W_\bullet:0\subsetneq W_1\subsetneq \cdots \subsetneq W_s=V_\alpha\ \text{with}\ \underline{dim}\ W_j/W_{j-1}=a_ji_j \},$$
 $$\widetilde{\mathscr F}_{\iii,\aaa}=\{(\underline{x}, W_\bullet)\mid \underline{x}\in E_{\alpha} \ \text{nilpotent}, W_\bullet\in  \mathscr F_{\iii,\aaa}\ \text{such that} \ \underline{x}(W_j)\subseteq W_{j-1}\}.$$
 Let $G_\alpha$ acts on $\widetilde{\mathscr F}_{\iii,\aaa}$  diagonally. The first projection $\pi_{\iii,\aaa}:\widetilde{\mathscr F}_{\iii,\aaa}\rightarrow E_{\alpha}$ is a $G_{\alpha}$-equivariant proper map, which yields $L_{\iii,\aaa}=({\pi_{\iii,\aaa}})_! (\K_{\widetilde{\mathscr F}_{\iii,\aaa}})[\text{dim}\widetilde{\mathscr F}_{\iii,\aaa}]$   a semisimple complex in $ D_{G_{\alpha}}(E_{\alpha})$. We set
 \vskip 1mm
  $\mathcal P_{\alpha}$:   the set of isomorphism classes of  simple perverse sheaves
  appearing, with possible shifts, in $L_{\iii,\aaa}$  for all $(\iii,\aaa)\vdash \alpha$,

 $\mathcal Q_{\alpha}$: the full subcategory of $ D_{G_{\alpha}}(E_{\alpha})$ whose objects are finite direct sums of shifts of the semisimple perverse sheaves coming from  $\mathcal P_{\alpha}$,

$\mathcal K_\alpha $: the Grothendieck group of  $\mathcal Q_{\alpha}$.
\vskip 1mm
Form $\mathcal K=\bigoplus_{\alpha\in \N[I]}\mathcal K_\alpha$. It was proved in \cite{Bozec2014b} that $\mathcal K$ has a geometrically defined (twisted) $\mathcal A$-bialgebra structure  that is isomorphic to the $_{\mathcal A} U^-$, associated to the symmetric Borcherds-Cartan matrix $A$ given by
 $$a_{ii}=2-2\ \#\{\text{loops on}\ i\},\ a_{ij}=-\#\{\text{arrows between} \ i \ \text{and}\ j \ \}\ \text{for}\ i\neq j.$$
 This isomorphism is given explicitly as follows
 $$F_i^{(a)}\leftrightarrow {\K}_{E_{ai}}\ \ \text{for}\ i\in I^+;\quad F_{ia}\leftrightarrow (\pi_{i,a})_!({\K}_{\{{0}\}})\ \ \text{for}\  i\in I^{\leq 0}, a> 0.$$
 We identify $\mathcal K$ and $_{\mathcal A} U^-$. So when $i\in I^{\leq 0}$, $F_{ia}$ are simple
 perverse sheaves supported on $\{0\}\subseteq E_{ai}$.

\begin{remark}\label{iso}
  For $i\in I^0, a>0$, we denote by $\{\mathcal O_{\lambda}\}_{\lambda\vdash a}$ the nilpotent orbits (labelled by partitions of $a$) in $E_{ai}^{nil}$ under the action of $GL_{ai}$. Then $F_{ia}=IC(\mathcal O_{(1^{a})})$ is the simple perverse sheaf associated to the closed  orbit $\{0\}$   in $E_{ai}^{nil}$.

  The power $F_{i}^a$ is   the Springer sheaf $\Spr_{GL_a}=\pi_! ({\K}_{\widetilde{\mathcal N}_a}[\text{dim}\widetilde{\mathcal N}_a])$, where $\mathcal N_a=E_{ai}^{nil}$ and $\pi:\widetilde{\mathcal N}_a\rightarrow \mathcal N_a$ is the Springer map. Therefore,  we can write  $F_{i}^a=\bigoplus_{\lambda\vdash a} IC(\mathcal O_{\lambda})\otimes V_{\lambda}$ for some nonzero vector spaces $V_{\lambda}$, and $F_{i}^a$ corresponding to the regular $\K[S_a]$-module under the Springer correspondence.
\end{remark}
 There is a  geometric pairing $\{ \ , \ \} \colon \mathcal K \times \mathcal K\rightarrow \Z (\!(q)\!)$ defined by the equivariant cohomology (see e.g. \cite[8.1.9]{Lus}) which is coincide with the one we define on the $_{\mathcal A} U^-$, especially, we have for all $i\in I, a>0$,
 \begin{equation}\{F_{ia}, F_{ia}\}=\sum_{j}\text{dim}\ H^{j}_{GL_{a}}{(\text{pt})}\ q^j=\prod_{k=1}^{a}\frac{1}{1-q^{2k}}.\end{equation}
Here, if $i\in I^+$, $F_{ia}=F_i^{(a)}$.

 \vskip 3mm
We also consider the subalgebra $\mathcal K^1$ of $\mathcal K$ defined as follows:
\vskip 1mm
  $\mathcal P_{\alpha}^1$:   the set of isomorphism classes of   simple perverse sheaves
  appearing in $L_{\iii,\aaa}$  for all $(\iii,\aaa)\vdash \alpha$ with each $a_j=1$,

 $\mathcal Q_{\alpha}^1$: the full subcategory of $ D_{G_{\alpha}}(E_{\alpha})$ generated by $\mathcal P_{\alpha}^1$,

 $\mathcal K^1_\alpha$: the Grothendieck group of  $\mathcal Q_{\alpha}^1$.
\vskip 1mm
 Then $\mathcal K^1=\bigoplus_{\alpha\in \N[I]}\mathcal K^1_\alpha$ is the subalgebra of $\mathcal K$ generated by  $ F_i$ for $i\in I^+\cup I^-$,  and $F_{ia}$ for $i\in I^{ 0}, a>0$.

  \vskip 3mm
 We mention here  a smaller subalgebra, which is known as the {\it quantum generalized Kac-Moody algebra} introduced in \cite{Kang}.
   For $(\iii,\aaa)\vdash \alpha$ with each $a_j=1$, we could identify $(\iii,\aaa)$ with $\iii$. We set $\iii^{\leq 0}=(i_{\ell_1},\dots, i_{\ell_p})$ to be the subsequence of imaginary indices in $\iii$,  and   see that
 $\iii^{\leq 0}\vdash \alpha^{\leq 0}=\sum_{i\in I^{\leq 0}}\alpha_i i$. Let $$ \mathscr F_{\iii^{\leq 0}}=\{\text{all flags}\  W_\bullet:0\subsetneq W_1\subsetneq \cdots \subsetneq W_p=V_{\alpha^{\leq 0}}\ \text{with}\ \underline{dim}\ W_j/W_{j-1}=i_{\ell_j} \},$$
  $$\widetilde{\mathscr F}_{\iii^{\leq 0}}=\{(\underline{x}, W_\bullet)\mid \underline{x}\in E_{\alpha} \ \text{nilpotent}, W_\bullet\in  \mathscr F_{\iii^{\leq 0}}\ \text{such that} \ \underline{x}(W_j)\subseteq W_{j-1}\oplus V_{\alpha^+}\},$$
  where $V_{\alpha^+}=\bigoplus_{i\in I^+}\K^{\alpha_i}$.  We have the commutative diagram of $G_{\alpha}$-equivariant map:
  $$\xymatrix{\widetilde{\mathscr F}_{\iii}\ar[r]^{\theta_{\iii}\ \ }\ar[dr]_{\pi_{\iii}}&\widetilde{\mathscr F}_{\iii^{\leq 0}}\ar[d]^{\pi'_{\iii}} \\ & E_{\alpha}}$$
  where $\theta_{\iii}: (\underline{x}, W_\bullet)\mapsto (\underline{x}, 0\subsetneq W_{i_{\ell_1}}^{\leq 0}\subsetneq \cdots \subsetneq W_{i_{\ell_p}}^{\leq 0}=V_{\alpha^{\leq 0}})$,  $\pi'_{\iii}$ be the  first projection, which is shown to be semismall in   \cite{KS06}.

Let $\tau_\iii$ be the set of  simple perverse sheaves appearing in $ ({\theta_{\iii}})_! ({\K}_{\widetilde{\mathscr F}_{\iii}}[\text{dim}\widetilde{\mathscr F}_{\iii}])$, and let
  $\mathcal P_{\alpha}^2=\bigsqcup_{\iii\vdash \alpha}\{({\pi'_{\iii}})_{!}(P)\mid P\in \tau_\iii \}$,  a set of semisimple perverse sheaves. We set
\vskip 1mm

 $\mathcal Q_{\alpha}^2$: the full subcategory of $ D_{G_{\alpha}}(E_{\alpha})$ generated by $\mathcal P_{\alpha}^2$,

  $\mathcal K^2_\alpha$: Grothendieck group of  $\mathcal Q_{\alpha}^2$.
\vskip 1mm
 Then $\mathcal K^2=\bigoplus_{\alpha\in \N[I]}\mathcal K^2_\alpha$ is the subalgebra of $\mathcal K$ generated by  $ F_i, i\in I$.

  \vskip 3mm

\begin{remark}
A KLR-categorification of $\mathcal K^2$  was studied in \cite{KOP2012} and \cite{TW2023}, particularly for an arbitrary Borcherds-Cartan datum in \cite{TW2023}. We expect that the construction in \cite{TW2023} corresponds to the `canonical' basis $\mathcal P^2$.
  The main goal of this work is to provide a KLR-categorification for $\mathcal K$ and $\mathcal K^1$.
\end{remark}

\vskip 6mm
\section{\textbf{Categorification of  ${U}^-$ }}


\vskip 3mm
\subsection{Generators $\mathtt b_{i\ell}$ for $i\in I^-$}\
\vskip 3mm

Given a Bocherds-Cartan datum $(I, A, \cdot)$,  let $U^-$ be the   associated  quantum Borcherds-Bozec algebra.

Let $^-$ be  the $\Q$-algebra involution of $U^-$ given by
$ \overline{F}_{i\ell}=F_{i\ell}$ for all $(i,\ell)\in I^\infty$ and $\overline{q}=q^{-1}
$. Let $^*$ be the $\Q(q)$-algebra anti-involution of $U^-$ given by
${F}_{i\ell}^*=F_{i\ell}$ for all $(i,\ell)\in I^\infty$.

\vskip 3mm
\begin{proposition}\label{Bozec}\cite{Bozec2014b,Bozec2014c}\ \ {\it Let $i\in I^-$.
There exists  a unique set of elements $\{\mathtt b_{i\ell}\mid \ell \geq 1\}$ in $U^-$, such that $\mathtt b_{i\ell}\in  U^-_{\ell i}$ and
\begin{itemize}
\item[(1)] $\mathtt b_{i\ell}-F_{i\ell}\in \Q (q) \left<F_{ik} \mid k<\ell \right>$,

\item[(2)] $\{\mathtt b_{i\ell},z\}=0$ for all $z\in \Q (q)  \left<F_{i1} ,\cdots,F_{i\,\ell-1}\right>$.
\end{itemize}
 For any $\mathbf c =(c_1,\dots,c_t)\in\N^t$,  we set $ \mathtt b_{i,\mathbf c}=\mathtt b_{ic_1}\cdots\mathtt b_{ic_t}$. These elements satisfy the following properties:
\begin{itemize}
\item[(i)] $\rho(\mathtt b_{i\ell})=\mathtt b_{i\ell}\otimes 1+ 1\otimes \mathtt b_{i\ell}$, $\overline{\mathtt b}_{i\ell}=\mathtt b_{i\ell}$ and ${\mathtt b}_{i\ell}^*={\mathtt b}_{i\ell}$,

\item[(ii)] $\{\mathtt b_{i\ell},\mathtt b_{i\ell}\}\equiv 1  \pmod{q}$,\\
 $\{\mathtt b_{i,\mathbf c},\mathtt b_{i,\mathbf c'}\}=0$ if $\mathbf c$ and $\mathbf c'$ determine  different partitions,

\item[(iii)]  the set $\{ \mathtt b_{i,\mathbf c}\mid \mathbf c\in\mathcal C_\ell\}$ forms a basis of  $U^-_{\ell i}$.
\end{itemize}}
\end{proposition}

\vskip 3mm
\begin{proposition}\label{Tong}\cite[Theorem 2.4]{FKKT2022}\  \ {\it We have an algebra automorphism $\Psi\colon U^-\rightarrow U^-$ given by
\begin{equation}\label{Psi}
\Psi(F_i)=F_i ,\ i\in I^+;\ \ \Psi(F_{i\ell})=F_{i\ell},\ i\in I^0, \ell\geq 1;\ \ \Psi(\mathtt b_{i\ell})=F_{i\ell} , \ i\in I^-, \ell\geq 1.
\end{equation}  }
\end{proposition}

\vskip 5mm

\subsection{Khovanov-Lauda-Rouquier algebras $R(\nu)$}\
\vskip 3mm

 Let the nondegenerate symmetric bilinear form $\{ \ , \ \}$ on $U^-$    take specific values for   $F_{i\ell}$ $((i,\ell)\in I^\infty)$ as follows:
$$
\{F_{i\ell},F_{i\ell}\}=\begin{cases} 1/(1-q_i^2)\qquad \text{for}\ i\in I^+,\\ 1/(1-q_i^2)(1-q_i^4)\cdots (1-q_i^{2\ell})\quad \text{for}\ i\in I^0, \ell\geq 1,\end{cases}$$
and for $i\in I^-$, the value $\{F_{i\ell}, F_{i\ell}\}$   leads to
$$\{\mathtt b_{i\ell}, \mathtt b_{i\ell}\}=1/(1-q_i^2)$$
for each $\ell\geq 1$. Note that this setting  satisfies the requirement $\{F_{i\ell}, F_{i\ell}\} \equiv 1  \pmod{q}$ for all $(i,\ell)\in I^\infty$.

\vskip 3mm

Define $\mathbb I=I^+\cup I^0\cup (I^{-}
\times \Z_{>0})$. 
Denote by $\mathcal X$  the set of all equivalence classes $[\ii]$ of
sequences in $\mathbb I$,  two sequences $\ii$ and $\jj$ are equivalent
 if they permute each others.

 Let $\nu=[\ii]=[(i_1,\ell_1)\dots (i_n,\ell_n)]\in \mathcal X$. We set $\ell(\nu)=n$ and $|\nu|=\ell_1{i_1}+\cdots+\ell_n{i_n}\in \N[I]$.
For $\nu=[(i_1,\ell_1)\dots (i_n,\ell_n)],\nu'=[(j_1,k_1)\dots (j_m,k_m)]$, we set
$$\nu+\nu'=\nu\cup \nu'=[(i_1,\ell_1)\dots (i_n,\ell_n),(j_1,k_1)\dots (j_m,k_m)].$$ If $\nu'=[(i_1,\ell_1)\dots (i_t,\ell_t)]$  is a part of $\nu=[(i_1,\ell_1)\dots (i_n,\ell_n)]$, then we set $$\nu-\nu'=\nu\backslash \nu'=[(i_{t+1},\ell_{t+1})\dots (i_n,\ell_n)].$$

We assign a graph $\Upsilon$ with vertices set $\mathbb I$ and an edge between $(i,\ell)\neq(j,k)$ when  $a_{ij}\neq 0$.
\vskip 2mm

\begin{definition}
Fix an $\nu\in \mathcal X$ of length $n$.
We define the Khovanov-Lauda-Rouquier algebra
$R(\nu)$ associated to a given  Borcherds-Cartan datum $(I,A,\cdot)$ to be the $\K$-algebra with the homogeneous generators given by diagrams
(see \cite{KL2009} for a detailed explanation of the braid-like planar diagrams):

$$\begin{aligned}&\ 1_{\ii}=\genO{(i_1,\ell_1)\ \ \ }{(i_k,\ell_k)}{ \ \ \ (i_n,\ell_n)} \quad \text{for} \ \ii=\ii_1\dots\ii_n=(i_1,\ell_1)\dots (i_n,\ell_n)\in \nu, \ \text{deg}(1_{\ii})=0,\\ &\\ & x_{k,\ii}=\genX{(i_1,\ell_1)\ \ }{(i_k,\ell_k)}{\ \ (i_n,\ell_n)} \quad \text{for} \ \ii\in \nu, 1\leq k\leq n,\ \text{deg}(x_{k,\ii})=2r_{i_k},\end{aligned}$$
\vskip 2mm
$$\tau_{k,\ii}=\genT{(i_1,\ell_1)\ \ }{(i_k,\ell_k)}{(i_{k+1},\ell_{k+1})}{(i_n,\ell_n)} \quad \text{for} \  \ii\in \nu, 1\leq k\leq n-1,\ \text{deg}(\tau_{k,\ii})=-\ell_k\ell_{k+1}r_{i_k}a_{i_ki_{k+1}}.$$
\vskip 3mm
\noindent subject to the following local relations:

\begin{align}\label{dcross}
 \dcross{(i,\ell)}{(j,k)} \  =  \  \begin{cases} \quad \quad \quad \quad \quad \  0 & \text{ if } i=j\in I^+, \\
                    \\ \quad  \quad \quad  \ \ \dcrossA{(i,\ell)}{(j,k)}  & \text{ if } a_{ij}=0, \\
                    \\
                    \   \Big( \dcrossLL{-\ell^2\frac{a_{ii}}{2}}{(i,\ell)}{(i,\ell)} \ + \ \dcrossRR{-\ell^2\frac{a_{ii}}{2}}{(i,\ell)}{(i,\ell)} \Big)^2 & \text{ if }i=j\in I^- \ \text{and} \  \ell=k,\\
                    \\
                    \   \dcrossL{-k\ell a_{ij}}{(i,\ell)}{(j,k)} \ + \ \dcrossR{-k\ell a_{ji}}{(i,\ell)}{(j,k)}  & \text{ if } (i,\ell)\ne (j,k) \ \text{and} \  a_{ij} \ne 0,
                  \end{cases}
\end{align}

\begin{equation}
\begin{aligned}
 \LU{{}}{i}{i} \   -  \ \RD{{}}{i}{i}\ =\ \dcrossAA{i}{i} \quad\quad\quad \LD{{}}{i}{i}   \ - \ \RU{{}}{i}{i}\ =\ \dcrossAA{i}{i}  \quad \text{ if } i \in I^+,
\end{aligned}
\end{equation}

\begin{equation}
\begin{aligned}
 \LU{{}}{(i,\ell)}{(j,k)} \   =  \ \RD{{}}{(i,\ell)}{(j,k)}  \quad\quad\quad \LD{{}}{(i,\ell)}{(j,k)}   \ =  \ \RU{{}}{(i,\ell)}{(j,k)}  \quad \text{otherwise},
\end{aligned}
\end{equation}

\begin{equation}
\quad \quad\quad\quad \quad\BraidLLL{i}{(j,\ell)}{i} \  -  \ \BraidRRR{i}{(j,\ell)}{i} \ =  {\sum_{c=0}^{-\ell a_{ij}-1}} \ \threeDotStrands{i}{(j,\ell)}{i} \ \ \text{ if } i\in I^+, i\ne j \ \text{and} \ a_{ij} \ne 0,
 \end{equation}
 \begin{equation}
\BraidL{(i,\ell)}{(j,k)}{(h,m)} \  = \ \BraidR{(i,\ell)}{(j,k)}{(h,m)}\quad\quad \text{otherwise}.
\end{equation}

\end{definition}

\vskip 3mm
For $\ii,\jj\in \nu$, we set $_{\jj}R(\nu)_{\ii}=1_{\jj}R(\nu)1_{\ii}$, then  $R(\nu)=\bigoplus_{\ii,\jj} {_{\jj}R(\nu)_{\ii}}$. Denote by $P_{\ii}=R(\nu)1_{\ii}$ (resp. $_{\jj}P=1_{\jj}R(\nu)$) the gr-projective left (resp. right) $R(\nu)$-module.

\vskip 3mm
 For $\ii\in \nu$, set $\mathscr{P}_{\ii}={\K}[x_1(\ii),\dots,x_n(\ii)]$ and form the $\K$-vector space $\mathscr{P}_{\nu}=\bigoplus_{\ii\in \nu}\mathscr{P}_{\ii}$. Each $\omega\in S_n$ acts on $\mathscr{P}_{\nu}$ by sending $x_a(\ii)$ to $x_{\omega (a)}(\omega(\ii))$.

\vskip 3mm

Choose an orientation for each edge of $\Upsilon$.
 We define an action of $R({\nu})$ on $\Pol_{\nu}$ as follows.

\vskip 1mm

\begin{itemize}

\item[(i)]  If $\ii\neq\kk$,  $_{\jj}R(\nu)_{\ii}$ acts on $\Pol_{\kk}$ by $0$.

\vskip 1mm

\item[(ii)]  For $f\in \Pol_{\ii}$,
$1_{\ii}\cdot f=f, \ \ x_{k,\ii}\cdot f=x_k{(\ii)}f.$

\vskip 1mm

\vskip 1mm

\item[(iii)] If $\ii_k=(i,\ell)$, $\ii_{k+1}=(j,t)$,
\begin{equation*}
\tau_{k,\ii}\cdot f=\begin{cases}\frac{f- s_kf}{x_k{(\ii)}-x_{k+1}(\ii)} & \text{if}\ i= j\in I^+,\\
{s_k}f & \text{if}\  a_{ij}=0 \ \text{or if}\ \ii_k\leftarrow \ii_{k+1}, \\
(x_k(\ii)^{-\ell^2\frac{a_{ii}}{2}}+x_{k+1}(\ii)^{-\ell^2\frac{a_{ii}}{2}})s_kf
& \text{if}\ i=j\in I^- \ \text{and}\ \ell=t,\\
 (x_k(s_k\ii)^{-\ell ta_{ji}}+x_{k+1}(s_k\ii)^{-\ell ta_{ij}})s_kf & \text{if}\ \ii_k\rightarrow \ii_{k+1}.
\end{cases}
\end{equation*}
\end{itemize}
It is  easy to check $\Pol_\nu$ is an $R(\nu)$-module with the action defined above.


\vskip 5mm

\subsection{Algebras $R(n(i,\ell))$ and their gr-irreducible modules}\

\vskip 3mm
Fix  $(i,\ell)\in \mathbb I$ and $n\geq 0$. Let $\nu=[n(i,\ell)]\in \mathcal X$ that has only one sequence $\ii=\underbrace{(i,\ell)\dots (i,\ell)}_{n}$.
The  algebra $R(n(i,\ell))$ is generated by $x_{1,\ii},\dots, x_{n,\ii}$ of degree $2r_i$
and $\tau_{1,\ii},\dots,\tau_{n-1,\ii}$ of degree $-\ell^2r_ia_{ii}$
subject to the following local relations:
 \vskip 1mm
 $$\dcrossI{i}{i}=0\quad\quad \ \LU{{}}{i}{i}   -   \RD{{}}{i}{i}\ =\LD{{}}{i}{i}   -  \RU{{}}{i}{i}\ =\ \dcrossAA{i}{i}\quad\quad \ \BraidLI{i}{i}{i}  =  \BraidRI{i}{i}{i} \quad\quad \ \text{ if } i \in I^+ .$$
 \vskip 2mm
 $$\dcrossII{(i,\ell)}{(i,\ell)}=\bigg( \dcrossLL{-\ell^2\frac{a_{ii}}{2}}{(i,\ell)}{(i,\ell)} \ + \ \dcrossRR{-\ell^2\frac{a_{ii}}{2}}{(i,\ell)}{(i,\ell)} \bigg)^2\quad\ \ \LU{{}}{(i,\ell)}{(i,\ell)}    =  \RD{{}}{(i,\ell)}{(i,\ell)}\quad \ \ \LD{{}}{(i,\ell)}{(i,\ell)}     =    \RU{{}}{(i,\ell)}{(i,\ell)}$$  $$ \BraidLII{(i,\ell)}{(i,\ell)}{(i,\ell)}  =  \BraidRII{(i,\ell)}{(i,\ell)}{(i,\ell)} \quad \text{ if } i \in I^- .$$
 \vskip 2mm
 $$\dcrossI{i}{i}\ = \  \dcrossAA{i}{i}\quad \quad \LU{{}}{i}{i}    =  \RD{{}}{i}{i}\quad \quad \LD{{}}{i}{i}   =  \RU{{}}{i}{i}\quad\quad  \BraidLI{i}{i}{i}  =  \BraidRI{i}{i}{i} \quad \quad \text{ if } i \in I^0 .$$

 \vskip 1mm

We will abbreviate $x_{k,\ii}$ (resp. $\tau_{k,\ii}$)  for
 $x_k$ (resp. $\tau_k$).
 In all cases, $R(n(i,\ell))$ has a basis
 $$ \{x_1^{r_1}\cdots x_n^{r_n}\tau_\omega \mid
 \omega \in S_n,r_1,\dots,r_n\geq 0\}.$$
Indeed, for instance, if $i\in I^0$ or $I^-$, we just need  to
consider the actions of these elements on $x_1^Nx_2^N\cdots x_n^{nN}$ for $N\gg 0$.

Therefore, we could identify the polynomial algebra $P_n=\K[x_1,\dots,x_n]$ with the subalgebra of $R(n(i,\ell))$ generated by $x_1,\dots,x_n$. Then the center of $R(n(i,\ell))$ is $Z_n$, the algebra of symmetric polynomials in $x_1,\dots,x_n$.

\vskip 3mm
 Let $i\in I^+$.  By the representation theory of the nil-Hecke algebras, $R(ni)$ has a unique gr-irreducible module $V(i^n)$ of graded dimension $[n]_i!$, which is isomorphic to $R(ni)\otimes_{P_n}\mathbf 1_n\left\{\frac{n(n-1)}{2} \cdot r_i\right\}$.
  Here, $\mathbf 1_n$ is the one-dimensional trivial module over $P_n$ on which  each $x_k$ acts by $0$. 
\vskip 3mm
Let $i\in I^-$, since $R(n(i,\ell))$ has only trivial idempotents, it has a unique gr-irreducible module $V((i,\ell)^n)$, which  is the one-dimensional trivial module with the gr-projective cover $R(n(i,\ell))$.
\vskip 3mm
Now let $i\in I^0$. Note that $R(ni)_0$ is just the symmetric group algebra $\K S_n$. Since $\text{char} \K=0$,  it is well known that  $\K S_n$ has $|\mathcal P_n|$ many irreducible modules that can be labelled by the partitions of $n$.

Let $V$ be an irreducible $\K S_n$-module (which is  also an indecomposable projective module). Then
$$\widetilde{V}:=R(ni)\otimes_{R(ni)_0}V \Big/ R(ni)_{> 0}\otimes_{R(ni)_0}V $$
is a gr-irreducible $R(ni)$-module. In other words, $\widetilde{V}$ is obtained from $V$ with the actions of $x_1,\dots,x_n$ by $0$.
Moreover, all gr-irreducible  $R(ni)$-modules can be obtained in this way.
When no confusion arises, we still write $V$ for $\widetilde{V}$.

\vskip 3mm

We have shown the following:
\begin{proposition}
{\it If $V_1,\dots,V_{|\mathcal P_n|}$ is a complete set of non-isomorphic classes of irreducible $\K S_n$-modules, then $V_1,\dots,V_{|\mathcal P_n|}$ is a complete set of non-isomorphic classes of gr-irreducible $R(ni)$-modules. In particular, the gr-Jacobson radical $J^{\text{gr}}(R(ni))=R(ni)_{>0}$.
}
\end{proposition}

\vskip 3mm
Let $V_{i,n}$ be the one-dimensional trivial module over $\K S_n$. Note that
$$V_{i,n}=\K S_n\cdot e_{i,n}=\K\cdot e_{i,n},$$
where $e_{i,n}=\frac{1}{n!}\sum_{\omega \in S_n} \omega$. If $r+t=n$, then the restriction  to  $\K S_r\otimes \K S_t$-modules gives
$$\Res^n_{r,t}V_{i,n}\cong V_{i,r}\otimes V_{i,t}=\K S_n\cdot e_{i,r}\otimes \K S_n\cdot e_{i,t}.$$
\vskip 3mm
The gr-projective cover of $R(ni)$-module $V_{i,n}$ is
$P_{i,n}=R(ni)e_{i,n}$, which has a basis
$$\{x_1^{r_1}\cdots x_n^{r_n}\cdot e_{i,n}\mid r_1,\dots,r_n\geq 0\},$$
 the restriction of $P_{i,n}$ to  $R(ri)\otimes R(ti)$-modules gives
\begin{equation}\label{P}\Res^n_{r,t}P_{i,n}\cong P_{i,r}\otimes P_{i,t}.\end{equation}

\vskip 3mm
Since $e_{i,n}R(ni)e_{i,n}$ is spanned by $\{f\cdot e_{i,n}\mid f\in Z_n\}$, we see that
\begin{equation}\label{P1}(P_{i,n},P_{i,n})=\gdim(e_{i,n}R(ni)e_{i,n})=\gdim \ Z_n=1/(1-q_i^2)(1-q_i^4)\cdots (1-q_i^{2n}).\end{equation}
Here, $( \ , \ )$ is the Khovanov-Lauda's form defined in (\ref{KL-form}).
\vskip 5mm
\subsection{Grothendieck groups $K_0(R)$ and $G_0(R)$}\

\vskip 3mm
Let $\ii,\jj\in \nu$. Using the polynomial representation $\mathscr P_{\nu}$ of $R(\nu)$, one can  obtain by a similar argument in \cite[Theorem 2.5]{KL2009} that $\mathscr P_{\nu}$ is a faithful  $R(\nu)$-module and  $_{\ii}R(\nu)_{\jj}$ has a basis $$\{x_{1,\ii}^{u_1}\cdots x_{n,\ii}^{u_n}\cdot\widehat\omega_\jj\mid u_1,\dots,u_n\in \N,\ \omega\in S_n\ \text{such that} \ \omega(\jj)=\ii\},$$ where $\widehat\omega_\jj\in {_{\ii}} R(\nu)_{\jj}$ is uniquely determined by a fixed  reduced expression of $\omega$.

\vskip 3mm

Assume $\nu$ contains a sequence $(i_1,\ell_1)^{m_1}\cdots(i_t,\ell_t)^{m_t}$ such that $(i_1,\ell_1),\dots,(i_t,\ell_t)$ are all distinct. Similar to \cite[Theorem 2.9]{KL2009}, the center $Z(R{(\nu)})$ of $R{(\nu)}$ can be  described as $$Z(R(\nu))\cong \bigotimes^t_{k=1}\K[z_{1},\dots,z_{m_k}]^{S_{m_k}},$$ the latter is a tensor product of symmetric polynomial algebras such that the generators in $\K[z_{1},\dots,z_{m_k}]$ are of degree $2r_{i_k}$. Moreover, $R(\nu)$ is a free $Z(R(\nu))$-module of rank $((m_1+\cdots+m_t)!)^2$. It is also a gr-free $Z(R(\nu))$-module of finite rank. So we have $$\gdim Z(R(\nu))=\prod ^t_{k=1}\left(\prod^{m_k}_{c=1}\frac{1}{1-q_{i_k}^{2c}}\right)$$ and $\gdim R(\nu)\in \Z[q,q^{-1}]\cdot \gdim Z(R(\nu))$.

\vskip 3mm
Denote by
$$\begin{aligned}
&  R(\nu)\text{-}\Mod\colon \text{the category of finitely generated graded}\  R(\nu) \text{-modules},\\
& R(\nu)\text{-}\fMod\colon\text{the category of finite-dimensional graded}\  R(\nu) \text{-modules},\\
& R(\nu)\text{-}\pMod\colon\text{the category of  projective objects in } R(\nu)\text{-}\Mod.
\end{aligned}$$

 Up to isomorphism and degree shifts, each $R(\nu)$ has only finitely many gr-irreducible modules, all of which are finite-dimensional and are irreducible $R(\nu)$-modules by forgetting the grading. Let $\B_{\nu}$ be the set of equivalence classes of gr-irreducible $R(\nu)$-modules. Choose one representative $S_b$ from
each equivalence class and denote by $P_b$ the gr-projective cover of $S_b$. The Grothendieck group $G_0(R(\nu))$ (resp. $K_0(R(\nu))$) of $R(\nu)$-$\fMod$ (resp. $R(\nu)$-$\pMod$) are free $\Z[q,q^{-1}]$-modules with $q[M]=[M\{1\}]$, and with a  basis $\{[S_b]\}_{b\in \B_{\nu}}$ (resp. $\{[P_b]\}_{b\in \B_{\nu}}$).
\vskip 3mm

  Let $R=\bigoplus_{\nu\in \mathcal X}R(\nu)$  and form
$$G_0(R)=\bigoplus_{\nu\in \mathcal X}G_0(R(\nu)),\ K_0(R)=\bigoplus_{\nu\in \mathcal X}K_0(R(\nu)).$$
The   $K_0(R)$ and $G_0(R)$ are equipped with  twisted bialgebras structure induced by the induction and restriction functors:
\begin{equation*}
\begin{aligned}
& \Ind^{\nu+\nu'}_{\nu,\nu'}\colon R(\nu)\otimes R(\nu')\text{-}\Mod\rightarrow R(\nu+\nu')\text{-}\Mod,\ M\mapsto R(\nu+\nu')1_{\nu,\nu'}\otimes_{R(\nu)\otimes R(\nu')}M,\\
&\Res^{\nu+\nu'}_{\nu,\nu'}\colon R(\nu+\nu')\text{-}\Mod\rightarrow R(\nu)\otimes R(\nu')\text{-}\Mod,\ N\mapsto 1_{\nu,\nu'}N,
\end{aligned}
\end{equation*}
where $1_{\nu,\nu'}=1_\nu\otimes 1_{\nu'}$. More precisely,  we set $|x|=|\nu|\in \N[I]$ for $x\in R(\nu)\text{-}\Mod$ and equip $K_0(R)\otimes K_0(R)$ (resp. $G_0(R)\otimes G_0(R)$)   with a twisted algebra structure via
$$(x_1\otimes x_2)(y_1\otimes y_2)=q^{-|x_2|\cdot |y_1|}x_1y_1\otimes x_2y_2,$$
then $\Res$ is a $\Z[q,q^{-1}]$-algebra homomorphism by Mackey's Theorem \cite[Proposition 2.18]{KL2009}.
\vskip 3mm


\vskip 3mm

The $K_0(R)$ and $G_0(R)$ are dual to each other with respect to the bilinear pairing
$( \ , \ )\colon K_0(R)\times G_0(R)\rightarrow \Z[q,q^{-1}]$ given by
\begin{equation}\label{KL-form}([P],[M])=\gdim (P^\psi\otimes_{R(\nu)}M)=\gdim (\HOM_{R(\nu)}(\overline{P},M)),\end{equation}
where  $\psi$ is the anti-involution of $R(\nu)$ obtained by flipping the diagrams about horizontal axis and it turns a left $R(\nu)$-module into right, $\overline{P}=\HOM(P,R(\nu))^\psi$. There is also a symmetric bilinear form $( \ , \ )\colon K_0(R)\times K_0(R)\rightarrow \Z (\!(q)\!)$ defined in the same way.

\vskip 3mm




\vskip 3mm
The $K_0(R)$ and $G_0(R)$ are $\Z[q,q^{-1}]$-modules dual to each other  with respect to the bilinear pairing
  defined in (\ref{KL-form}). By (\ref{P1}) and  \cite[Proposition 3.3]{KL2009}, the symmetric bilinear form $( \ , \ )\colon K_0(R)\times K_0(R)\rightarrow \Z (\!(q)\!)$  satisfies
\begin{itemize}
\item[(1)] $([M],[N])=0$  \ if $M\in R(\nu)\text{-}\Mod$, $N\in R(\mu)\text{-}\Mod$ with $\nu\neq \mu$.
\item[(2)] $(1,1) = 1$, where $1=\K$ as a module over $R(0)=\K$.
\item[(3)] $([P_{i}], [P_{i}]) =1/(1-q_i^2)$ \ \ for $i\in I^+$ and $P_i=R(i)1_i$;\\
$([P_{(i,\ell)}], [P_{(i,\ell)}]) =1/(1-q_i^2)$ \ \ for $i\in I^-$ and $P_{(i,\ell)}=R((i,\ell))1_{(i,\ell)}$;\\
$([P_{i,\ell}], [P_{i,\ell}]) =1/(1-q_i^2)(1-q_i^4)\cdots (1-q_i^{2\ell})$ \ \ for $i\in I^0$ and $P_{i,\ell}=R(\ell i)e_{i,\ell}$.
\item[(4)] $(x, yz) = (\Res(x), y \otimes z)$  \ for $x,y,z
\in K_0(R)$.
\end{itemize}

\vskip 5mm
\subsection{Quantum Serre relations}\

\vskip 3mm

Let $\ii$ be a sequence with {\it divided powers}:
$$\ii=(j_1,a_1)\dots (j_{p_0},a_{p_0})i_1^{(m_1)}(k_1,b_1)\dots (k_{p_1},b_{p_1})i_2^{(m_2)}\dots i_t^{(m_t)} (h_1,c_1)\dots (h_{p_t},c_{p_t}), $$
where $i_1,\dots,i_t\in I^+$ and the others  belong to $\mathbb I$.
\vskip 3mm

For such an $\ii$, we assign the following idempotent.
\begin{equation*}
\begin{aligned}
1_{\ii}= & 1_{(j_1,a_1)\dots (j_{p_0},a_{p_0})}\otimes e_{i_1,m_1}\otimes 1_{(k_1,b_1)\dots (k_{p_1},b_{p_1})}\otimes e_{i_2,m_2} \\
& \otimes \cdots \otimes e_{i_t,m_t}\otimes  1_{(h_1,c_1)\dots (h_{p_t},c_{p_t})},\end{aligned}
\end{equation*}
where $e_{i,m}=x_1^{m-1}x_2^{m-2}\cdots x_{m-1}\tau_{w_0}$
with $w_{0}$ being the longest element in $S_{n}$.

\vskip 3mm

Set
\begin{equation*}
\begin{aligned}
\langle\ii \rangle & =\sum_{k=1}^t\frac{m_k(m_k-1)}{2}r_{i_k}, \\
 _\ii P & =1_\ii R(\nu)\{-\langle\ii \rangle\},\  P_\ii = R(\nu)\psi(1_\ii)\{-\langle\ii \rangle\}.
\end{aligned}
\end{equation*}
In particular, for $i\in I^+$ and $n\geq 0$,
 \begin{equation}\label{e}
 P_{i^{(n)}}= R(ni)\psi{(e_{i,n})}\left\{-\frac{n(n-1)}{2} \cdot r_i\right\} \cong R(ni)e_{i,n}\left\{ \frac{n(n-1)}{2} \cdot r_i \right\} .
 \end{equation}

\vskip 3mm
\begin{proposition}\label{Serre}
{\it Suppose $i\in I^+$, $j\in I$, $i\neq j$ and let $n\in \Z_{>0}$ and $m=1-na_{ij}$.
 Then we have isomorphisms of graded left $R(\nu)$-modules
\begin{align*}
&\bigoplus^{\lfloor \frac{m}{2} \rfloor}_{c=0}P_{i^{(2c)}j^ni^{(m-2c)}}\cong\bigoplus^{\lfloor \frac{m-1}{2} \rfloor}_{c=0}P_{i^{(2c+1)}j^ni^{(m-2c-1)}}\quad\text{if}\ j\in I^+,\\
&\bigoplus^{\lfloor \frac{m}{2} \rfloor}_{c=0}P_{i^{(2c)}(j,n)i^{(m-2c)}}\cong\bigoplus^{\lfloor \frac{m-1}{2} \rfloor}_{c=0}P_{i^{(2c+1)}(j,n)i^{(m-2c-1)}}\quad\text{if}\ j\in I^-,\\
& \bigoplus^{\lfloor \frac{m}{2} \rfloor}_{c=0}R(i^mj^n)\psi(e_{i,2c}\otimes e_{j,n}\otimes e_{i,m-2c})\langle i^{(2c)}i^{(m-2c)} \rangle\\
&\phantom{aaa} \cong\bigoplus^{\lfloor \frac{m-1}{2} \rfloor}_{c=0}R(i^mj^n)\psi(e_{i,2c+1}\otimes e_{j,n}\otimes e_{i,m-2c-1})\langle i^{(2c+1)}i^{(m-2c-1)}\rangle\quad\text{if}\ j\in I^0.
\end{align*}
Moreover,  if $a_{ij}=0$, then
$$\begin{aligned}& P_{(i,\ell)(j,k)}\cong P_{(i,\ell)(j,k)}\quad\quad\text{if}\ (i,\ell),(j,k)\in \mathbb I,\\& R(i^\ell(j,k))\cdot e_{i,\ell}\otimes 1_{(j,k)}\cong R(i^\ell(j,k))\cdot 1_{(j,k)}\otimes e_{i,\ell}\quad\quad\text{if}\ i\in I^0, (j,k)\in \mathbb I,\\& R(i^\ell j^k)\cdot e_{i,\ell}\otimes e_{j,k}\cong R(i^\ell j^k)\cdot e_{j,k}\otimes e_{i,\ell}\quad\quad \text{if}\ i,j \in I^0.
\end{aligned}$$}
\begin{proof}
The proof is the same as the `Box' calculations in \cite{KL2011}. We only explain the last isomorphism. Let $i,j\in I^0$ with $a_{ij}=0$. Note that
$$\Down\ \ \ \Up.$$
The right multiplication by
$$\DDown$$
is a map from $R(i^\ell j^k)\cdot e_{i,\ell}\otimes e_{j,k}$ to $R(i^\ell j^k)\cdot e_{j,k}\otimes e_{i,\ell}$, which has the obvious inverse by flipping this diagram.
\end{proof}
\end{proposition}

\vskip 3mm
Let $K_0(R)_{\Q(q)}=\Q(q)\otimes_{\Z[q,q^{-1}]}K_0(R)$. By (\ref{P}), Proposition \ref{Tong} and Proposition \ref{Serre}, we have a well-defined bialgebra homomorphism
$$\begin{aligned}
& \Gamma_{\Q(q)}\colon U^-\rightarrow K_0(R)_{\Q(q)}\\
&\phantom{\Gamma_{\Q(q)}:} F_i\mapsto [P_i]\   \qquad\quad \text{for}\ i\in I^+\\
& \phantom{\Gamma_{\Q(q)}:}F_{i\ell}\mapsto [P_{i,\ell}]\qquad\ \ \text{for}\ i\in I^0, \ell\geq 1\\
& \phantom{\Gamma_{\Q(q)}:}\mathtt b_{i\ell}\mapsto [P_{(i,\ell)}]\qquad \text{for}\ i\in I^-, \ell\geq 1\\
\end{aligned}$$
Now, the bilinear form $\{ \ , \ \}$ on $U^-$ and the Khovanov-Lauda's form  $( \ , \ )$ on $K_0(R)_{\Q(q)}$ coincide under the map $\Gamma_{\Q(q)}$, that is
$$(\Gamma_{\Q(q)}(x),\Gamma_{\Q(q)}(y))=\{x,y\} \ \ \text{for} \ x,y\in U^-.$$
Thus $\Gamma_{\Q(q)}$ is injective by the non-degeneracy of $\{ \ , \ \}$. Moreover, we have
$$ \Gamma_{\Q(q)}(\overline{x})=\overline{\Gamma_{\Q(q)}(x)}.$$

Let $_{\mathcal A} U^-_1$ be the $\mathcal A$-subalgebra of $U^-$ generated by  $ F_i^{(n)}$ for $i\in I^+$, $ F_{i\ell}$ for $i\in I^0,\ell\geq 1$ and  $\mathtt b_{i\ell}$ for $i\in I^-,\ell\geq 1$. Then $\Gamma_{\Q(q)}$ induces an injective $\Z[q,q^{-1}]$-bialgebra homomorphism
$\Gamma\colon{_{\mathcal A}}U^-_1 \rightarrow K_0(R)$.

\vskip 5mm
\subsection{Surjectivity of $\Gamma_{\Q(q)}$ and $\Gamma$}\

\vskip 3mm

Let $\nu\in \mathcal X$. Define $\underline{\nu}$ to be the set of sequence $\ii$ of type $\nu$ with `parameters' for $i\in I^0$. Such a sequence  is of the form
 $$\ii=(j_1,a_1)\dots (j_{p_0},a_{p_0})\bm{(i_1,n_1)}(k_1,b_1)\dots (k_{p_1},b_{p_1})\bm{(i_2,n_2)}\ \dots\ \bm{(i_t,n_t)}(h_1,c_1)\dots (h_{p_t},c_{p_t}), $$
with $(i_1,n_1),\dots,(i_t,n_t)\in I^0\times \Z_{>0}$ and such that the expended sequence
 $$(j_1,a_1)\dots (j_{p_0},a_{p_0})\underbrace{i_1\dots i_1}_{n_1}(k_1,b_1)\dots (k_{p_1},b_{p_1})\underbrace{i_2\dots i_2}_{n_2}\ \dots\ \underbrace{i_t\dots i_t}_{n_t} (h_1,c_1)\dots (h_{p_t},c_{p_t})$$
belongs to $\nu$. For each $\ii\in \underline{\nu}$,
we assign the following idempotent of $R(\nu)$
$$1_{\ii}=1_{(j_1,a_1)\dots (j_{p_0},a_{p_0})}\otimes e_{i_1,n_1}\otimes 1_{(k_1,b_1)\dots (k_{p_1},b_{p_1})}\otimes e_{i_2,n_2}\otimes \cdots \otimes e_{i_t,n_t}\otimes 1_{(h_1,c_1)\dots (h_{p_t},c_{p_t})}.$$


 \vskip 3mm
 We define the character of $M\in R(\nu)\text{-}\fMod$ as $$\Ch M=\sum_{\ii\in \underline{\nu}}\gdim (1_{\ii}M)\ii\ \in \Z[q,q^{-1}]\underline{\nu}.$$

\vskip 3mm
Each $\ii\in \underline{\nu}$ determines a monomial $\Theta_{\ii}$ in $U^-$ under the correspondence
$$i\mapsto F_i\ (i\in I^+), \quad (i,\ell)\mapsto F_{i \ell}\ (i\in I^0), \quad (i,\ell)\mapsto \mathtt b_{i \ell}\ (i\in I^-).$$
 Let $U^-_{\nu}$ be the $\Q(q)$-subspace of $U^-$ spanned by $\Theta_{\ii}$ for all $\ii\in \underline{\nu}$. Combining with $\Gamma_{\Q(q)}$, we obtain a $\Q(q)$-linear map
 $$\Q(q)\underline{\nu} \longrightarrow U^-_{\nu}\stackrel{\Gamma_{\Q(q)}}{\longrightarrow}K_0(R(\nu))_{\Q(q)},$$
 which has the dual map
$$G_0(R(\nu))_{\Q(q)}\stackrel{\Ch}{\longrightarrow}\Q(q)\underline{\nu}.$$
We next show that the character map $\Ch$ is injective.

\vskip 3mm
Let $i\in I^0$ and $\nu=ni$. In this case,
$$\underline{\nu}=\{(i,\ell_1)\dots(i,\ell_s)\mid(\ell_1,\dots,\ell_s)\in \mathcal C_n\} \ \ \text{and} \ \ U^-_{\nu}=U^-_{-n\alpha_i}. $$
Since $\Gamma_{\Q(q)}\colon U^-_{-n\alpha_i}\rightarrow K_0(R(\nu))_{\Q(q)}$ is injective and they have the same dimension $|\mathcal P_n|$, we see that $\Gamma_{\Q(q)}$ is an isomorphism and therefore $\Ch\colon G_0(R(ni))_{\Q(q)}\rightarrow\Q(q)\underline{\nu}$ is injective.


\vskip 3mm

\begin{lemma}\label{LLL1}
{\it Let $i\in I^0$. The characters of all non-isomorphic gr-irreducible $R(ni)$-modules are $\Q(q)$-linear independent.
}
\end{lemma}

\vskip 3mm
\begin{example} The irreducible $\K S_3$-modules and their characters are given by
$$\begin{aligned}
& S^{(3)}=P_{i,3}=\K S_3\cdot e_{i,3},\quad \Ch \ S^{(3)}=e_{i,3}+e_{i,1}\otimes e_{i,2}+e_{i,2}\otimes e_{i,1}+e_{i,1}\otimes e_{i,1}\otimes e_{i,1}\\
& S^{(21)}=\K S_3\cdot 1/3(1+s_1-s_2s_1-s_1s_2s_1), \quad \Ch\  S^{(21)}=e_{i,1}\otimes e_{i,2}+e_{i,2}\otimes e_{i,1}+2e_{i,1}\otimes e_{i,1}\otimes e_{i,1}\\
& S^{(111)}=\K S_3\cdot 1/6(1-s_1-s_2+s_1s_2+s_2s_1-s_1s_2s_1),\quad \Ch \ S^{(111)}=e_{i,1}\otimes e_{i,1}\otimes e_{i,1}\end{aligned}.$$
Here $S^{(3)},S^{(21)},S^{(111)}$ are the Specht modules.
\end{example}
\vskip 3mm

Let  $(i,\ell)\in \mathbb I$ and $n\geq 0$. Define  a functor
\begin{equation*}
\begin{aligned}
\Delta_{(i,\ell)^n}\colon R(\nu)\text{-}\Mod &\rightarrow R(\nu\backslash n(i,\ell))\otimes R(n(i,\ell))\text{-}\Mod \\
M & \longmapsto (1_{\nu\backslash n(i,\ell)}\otimes 1_{n(i,\ell)})M.
\end{aligned}
\end{equation*}
For each $M\in R(\nu)$-$\fMod$, we define
$$\varepsilon_{(i,\ell)}M=\Max\{n\geq 0\mid \Delta_{(i,\ell)^n}M\neq 0\}.$$

\vskip 3mm


The following lemma  can be proved by the same manner in \cite[Section 3.2]{KL2009}.

\vskip 3mm
\begin{lemma}\label{LLL2}
{\it Let $(i,\ell)\in \mathbb I$ and  $M\in R(\nu)$-$\fMod$ be a gr-irreducible module with $\varepsilon_{(i,\ell)}M=n$. Then $\Delta_{(i,\ell)^n}M$ is isomorphic to $K\otimes V$ for some gr-irreducible $K\in R(\nu\backslash n(i,l))$-$\fMod$ with $\varepsilon_{(i,\ell)}K=0$ and some gr-irreducible $V\in R(n(i,\ell))$-$\fMod$. Moreover, we have $$M\cong \hd\ \Ind_{\nu\backslash n(i,\ell),n(i,\ell)} K\otimes V.$$
}
\end{lemma}

\vskip 3mm
Recall that  a gr-irreducible $R(n(i,\ell))$-module $V$
 has one of the following forms:

\begin{itemize}

\item[(i)] if $i \in I^{+}$, then
 $V=V(i^n)$, the unique gr-irreducible module of the nil-Hecke algebra,

\item[(ii)] if  $i\in I^-$, then $V=V((i,\ell)^n)$ is the one-dimensional trivial module,

\item[(iii)] if $i \in I^{0}$, then $V$ has $|\mathcal P_n|$ many choices.

\end{itemize}


\vskip 3mm

\begin{theorem}
{\it The map $\Ch\colon G_0(R(\nu))_{\Q(q)}\rightarrow\Q(q)\underline{\nu}$ is injective.}
\begin{proof}
We show that  the characters of elements in $\B_{\nu}$ are linearly independent over $\Q(q)$ by induction on $\ell(\nu)$. The case of $\ell(\nu)=0$ is trivial. Assume for $\ell(\nu)<n$, our assertion is true. Now, suppose $\ell(\nu)=n$ and we are given a non-trivial linear composition
\begin{equation}\label{Ch}
\sum_{M}c_M\Ch M=0
\end{equation}
for some $M\in \B_{\nu}$ and some $c_M\in\Q(q)$. Choose  $(i,\ell)\in \mathbb I$. We prove by a downward induction on $k=n,\dots,1$ that $c_M=0$ for all $M$ with $\varepsilon_{(i,\ell)}M=k$.

If $k=n$ and $M\in \B_{\nu}$ with $\varepsilon_{(i,\ell)}M=n$, then $\nu=n(i,\ell)$ and $M$ is  a gr-irreducible $R(n(i,\ell))$-module. When $i\in I^+\sqcup I^-$, our assertion is trivial. When $i\in I^0$, it follows from Lemma \ref{LLL1}.

 Assume for $1\leq k<n$, we have $c_L=0$ for all $L$ with $\varepsilon_{(i,\ell)} L>k$. Taking out the terms with ${(i,\ell)}^k$-tail in the rest of (\ref{Ch}), we obtain
\begin{equation}\label{Ch1}\sum_{M\colon\varepsilon_{(i,\ell)}M=k}c_M\Ch (\Delta_{(i,\ell)^k}M)=0.\end{equation}
By Lemma \ref{LLL2}, we can assume $\Delta_{(i,\ell)^k}M\cong K_M\otimes V_M$ for  gr-irreducible $K_M\in R(\nu\backslash k(i,\ell))$-$\fMod$  with $\varepsilon_{(i,\ell)}K_M=0$ and gr-irreducible $V\in R(k(i,\ell))$-$\fMod$, then (\ref{Ch1}) becomes $$\sum_{M\colon\varepsilon_{(i,\ell)}M=k}c_M\Ch K_M\otimes \Ch V_M=0.$$
Note that if $[M]\neq [M']$ in $\B_{\nu}$, then we have $[K_M]\neq [K_{M'}]$ or $[V_M]\neq [V_{M'}]$.

By the inductive hypothesis, $\Ch K$ $(K\in \B_{\nu\backslash k(i,\ell)})$ are linearly independent, and by  Lemma \ref{LLL1} if $i\in I^0$, $\Ch V$ $(V\in \B_{ ki})$ are linearly independent. It follows that  $c_M=0$ for all $M$ with $\varepsilon_{(i,\ell)}M=k$.
Since each gr-irreducible $R(\nu)$-modules $M$ has $\varepsilon_{(i,\ell)}M>0$ for at least one ${(i,\ell)}\in \mathbb I$, the theorem has been proved.
\end{proof}
\end{theorem}
\vskip 3mm
\begin{remark}
We see from the  proof that the map $\Ch\colon G_0(R(\nu))\rightarrow\Z[q,q^{-1}]\underline{\nu}
$ is also injective. If we set $\text{ch}={\Ch |}_{q=1}$, then by a similar argument, the ungraded  characters  of elements in $\B_{\nu}$ are linearly independent over $\Z$.
\end{remark}
\vskip 3mm

\begin{corollary}{\it
$\Gamma_{\Q(q)}\colon  U^-\rightarrow K_0(R)_{\Q(q)}$ is an isomorphism.
}\end{corollary}


\vskip 3mm
We next consider the surjectivity of $\Gamma\colon {_{\mathcal A}}U^-_1\rightarrow K_0(R)$.
\vskip 3mm
For $\lambda\vdash n$, let $S^\lambda$ be the Specht module  corresponding to the Young diagram of shape $\lambda$. For $\mu\prec \lambda$, we denote by $S^{\lambda/\mu}$ the shew representation of $\K S_n$ corresponding to the skew diagram $\lambda/\mu$.
\vskip 3mm
\begin{lemma}\cite[Proposition 3.5.5]{CST2010}\label{LLL3}
{\it Let $\lambda\vdash n$ be a partition and $ (b_1,\dots,b_\ell) \vDash n$ be a composition of $n$. Then
$$\Res^n_{b_1,\dots,b_\ell}=\bigoplus(S^{\lambda^{(1)}}\otimes S^{\lambda^{(2)}/\lambda^{(1)}}\otimes \cdots\otimes S^{\lambda/\lambda^{(\ell-1)}}),$$
where the sum runs over all sequences $\lambda^{(1)}\prec \lambda^{(2)}\prec \cdots\prec \lambda^{(\ell)}=\lambda$ such that  $|\lambda^{(j)}/\lambda^{(j-1)}|=b_j$ for all $j=1,\dots,\ell$.
}
\end{lemma}

\vskip 3mm
\begin{lemma}\cite[Proposition 3.5.12]{CST2010}\label{LLL4}
{\it Assume $|\lambda/\mu|=k$. The multiplicity of the trivial representation $S^{(k)}$ in $S^{\lambda/\mu}$ is $1$ if $\lambda/\mu$ is totally disconnected, $0$ otherwise.
}
\end{lemma}
\vskip 3mm
\begin{proposition}\label{PPP1}
{\it Let $\lambda\vdash n$. Let  $\mathbf c \vDash n$ be a composition of $n$, which determines a partition $\lambda_{\mathbf c}\vdash n$. Then
$$\dim_\K( e_{i,\mathbf c}\cdot S^\lambda)=\begin{cases}1\quad \text{if}\ \lambda_{\mathbf c}=\lambda,\\ 0 \quad \text{if}\ \lambda_{\mathbf c}>\lambda,\end{cases}$$
where $>$ is the lexicographical order of partitions.
}\begin{proof}
The proposition follows from Lemma \ref{LLL3} and Lemma \ref{LLL4}, using the fact that the Kostka number
$$K_{\lambda, \mathbf c}=K_{\lambda, \lambda_{\mathbf c}}=\begin{cases}1\quad \text{if}\ \lambda_{\mathbf c}=\lambda,\\ 0 \quad \text{if}\ \lambda_{\mathbf c}>\lambda.\end{cases}$$
\end{proof}
\end{proposition}

Note that Lemma \ref{LLL1} can be derived directly from the proposition above.
\vskip 3mm
For $i\in I^0$ and $\lambda=(c_1,\dots,c_r)\vdash n$, we set $$P_{i,\lambda}=P_{i,c_1}\cdots P_{i,c_r}=R(ni)\cdot e_{i,\lambda}.$$
Then $(P_{i,{\lambda}}, S^{\mu}) 
=\text{dim}_\K(e_{i,\lambda}\cdot S^\mu)$ and according to Proposition \ref{PPP1}, the matrix $$\left\{(P_{i,{\lambda}}, S^{\mu})\right\}_{\lambda,\mu \vdash n}$$
 unitriangular.
It follows that
each $[P]\in K_0(R(ni))$ can be written as a $\Z[q,q^{-1}]$-linear combination of $[P_{i,\lambda}]$ for $\lambda \vdash n$.

\vskip 3mm
More generally, we write the set $\mathbb I$ as
$$\mathbb I=\{(i_1,\ell_1),(i_2,\ell_2),\dots,(i_k,\ell_k),\dots\}.$$
For a gr-irreducible $R(\nu)$-module $M$, we let $c_1=\varepsilon_{(i_1,\ell_1)}M$ and assume  $\Delta_{(i_1,\ell_1)^{c_1}}M = M_1\otimes V_M$ for  gr-irreducible $M_1$  with $\varepsilon_{(i_1,\ell_1)}M_1=0$ and gr-irreducible $V_M\in R(c_1(i_1,\ell_1))$-$\Mod$. If $i_1\in I^0$, then $V_M=S^{\lambda^{(1)}}$ for some $\lambda^{(1)}\vdash c_1$. So we get a pair $(c_1,\lambda^{(1)})$,  where we set $\lambda^{(1)}=0$ when $i_1\notin I^0$. Inductively,  $c_k=\varepsilon_{(i_{k},\ell_k)}M_{k-1}$ and $\Delta_{(i_k,\ell_k)^{c_k}}M_{k-1} = M_{k}\otimes V_{M_{k-1}}$, and we obtain $(c_k,\lambda^{(k)})$.
If we do not get $M_k=0$ after $\mathbb I$ exhausted, we can continue
the  above process from $(i_1,\ell_1)$.
Therefore, each $b\in \B_{\nu}$
is assigned by a sequence
$$W_b=(c_1,\lambda^{(1)})(c_2,\lambda^{(2)})\cdots (c_k,\lambda^{(k)})\cdots,$$
and we see from Lemma \ref{LLL2} that $b$ is uniquely determined by $W_b$.

\vskip 2mm

 Set
$$P_{(c_k,\lambda^{(k)})}=\begin{cases} P_{{i_k}^{(c_k)}}\qquad \text{if}\ i_k \in I^+,\\ P_{{(i_k,\ell_k)}^{c_k}}\quad \text{if}\ i_k\in I^-,\\ P_{{i_k}, \lambda^{(k)}} \quad\ \ \text{if}\ i_k\in I^0,\end{cases}$$
and
$P_{W_b}=\cdots P_{(c_k,\lambda^{(k)})} \cdots P_{(c_2,\lambda^{(2)})} P_{(c_1,\lambda^{(1)})}$.

\vskip 3mm

Let  $b,b'\in \B_{\nu}$ with $W_b=(c_1,\lambda^{(1)})(c_2,\lambda^{(2)})\cdots$ and $W_{b'}=(d_1,\mu^{(1)})(d_2,\mu^{(2)})\cdots$. We denote $W_b>W_b'$  if for some $t$, $(c_1,\lambda^{(1)})=(d_1,\mu^{(1)}),\dots,(c_{t-1},\lambda^{(t-1)})=(d_{t-1},\mu^{(t-1)})$ but $$c_t>d_t\ \text{or}\ c_t=d_t, \lambda^{(t)}>\mu^{(t)}.
$$
\vskip 3mm

\begin{proposition}
{\it $\HOM(P_{{W_b}}, S_{b'})=0$ if $b>b'$ and $\HOM(P_{{W_b}}, S_b)\cong\K$, up to a degree shift.
}
\begin{proof}
For $i\in I^+$,  we have $\HOM(P_{i^{(n)}},V(i^n))\cong \K$ since $P_{i^{(n)}}$ is the graded projective cover of $V(i^n)$. For $i\in I^-$, $\HOM(R({n(i,\ell)}),V((i,\ell)^n))\cong \K$ as graded vector spaces.
The results follows immediately from the Frobenius reciprocity and Proposition \ref{PPP1}, which deals with the case  $i\in I^0$.
\end{proof}
\end{proposition}

\vskip 3mm

By proposition above, each $[P]\in K_0(R(\nu))$ can be written as a $\Z[q,q^{-1}]$-linear combination of $[P_{{W_b}}]$ for $b\in \B_{\nu}$. Therefore, $\Gamma$ is surjective.
\vskip 3mm
\begin{theorem}{\it
 $\Gamma\colon {_\mathcal A} U^-_1 \rightarrow K_0(R)$ is an isomorphism.
}\end{theorem}
\vskip 3mm
For $M\in R(\nu)$-$\fMod$, let $M^\divideontimes=\HOM_\K(M,\K)^\psi$ be the dual module in $R(\nu)$-$\fMod$ with the action given by
$$(zf)(m):=f(\psi(z)m)\ \text{for}\ z\in R(\nu), f\in \HOM_\K(M,\K), m\in M.$$
 As proved in \cite[Section 3.2]{KL2009}, for each gr-irreducible $R(\nu)$-module $S$, there is a unique $r\in\Z$ such that $(L\{ r \})^\divideontimes \cong L\{r\}$, and the graded projective cover of $L\{ r \}$ is stable under the bar-involution $^-$.
\vskip 3mm
Recall that $_{\mathcal A} U^-$ is the $\mathcal A$-subalgebra of $U^-$ generated by  $ F_i^{(n)}$ for $i\in I^+$, $ F_{i\ell}$ for $i\in I^{\leq 0},\ell\geq 1$.
 Combined with $\Psi$ given in (\ref{Psi}), we  obtain a  $\Z[q,q^{-1}]$-algebra isomorphism
$$\Phi=\Psi\Gamma^{-1}\colon K_0(R)\xrightarrow{\sim}{_\mathcal A} U^-.$$

\begin{conjecture}
Under $\Phi$, the bar-invariant indecomposable projective modules of $K_0(R)$ coincides with the canonical basis   $\bigsqcup_{\alpha\in\N[I]}\mathcal P_{\alpha}$ of $_{\mathcal A}U^-$.
\end{conjecture}

We show this conjecture for the quiver with one vertex in the following section.
\vskip 5mm
\subsection{One vertex cases}\
\vskip 3mm
Assume   $I=I^{\leq 0}=\{i\}$ and ${_\mathcal A}U^-$ be the quantum Borcherds-Bozec associated to $I$. When $i\in I^-$,   the canonical bases are  the monomials in $F_{i\ell}$'s (see \cite{Lus93}). We see  that $\Phi$ maps the  self-dual indecomposable projective modules  of $K_0(R)$, which are the monomials in $[P_{(i,\ell)}]$'s, to the canonical bases.
\vskip 3mm
We next assume $i\in I^0$. Following the notations in  Remark \ref{iso}, it is already known  that  the canonical bases of ${_\mathcal A}U^-$  is  the $ \{IC(\mathcal O_{\lambda})\}_{\lambda\vdash n}$.  In the following lemma, we denote by $K_0(S_n)$ the Grothendieck group of the finite dimensional $\K[S_n]$-modules, then $\bigoplus_n K_0(S_n)$ is a $\Z$-algebra with the multiplication induced by the induction of modules.

\begin{lemma}\label{springer}
{\it If we have an isomorphism $\Omega: \mathcal A\otimes_{\Z}(\bigoplus_n K_0(S_n))\rightarrow {_\mathcal A} U^-$ of $\mathcal A$-algebras, which sends the trivial representation $S^{(n)}$     to $F_{in}=IC(\mathcal O_{(1^n)})$ or $IC(\mathcal O_{(n)})$, then it maps the irreducible $\K[S_n]$-modules to the canonical bases of ${_\mathcal A} U^-$. More precisely:
\begin{itemize}

\item[(i)] If $\Omega(S^{(n)})=F_{in}$, then $\Omega(S^\lambda)=IC(\mathcal O_{\widetilde{\lambda}})$ for all $\lambda\vdash n$, where   $\widetilde{\lambda}$ is the transpose of $\lambda$,

\item[(ii)] If  $\Omega(S^{(n)})=IC(\mathcal O_{(n)})$, then $\Omega(S^\lambda)=IC(\mathcal O_{\lambda})$ for all $\lambda\vdash n$.
\end{itemize}
}
\begin{proof}
Let $\Lambda$ be the ring of symmetric functions. By \cite{Lus81} (see also \cite[Example 3.10]{S09}), there is an $\mathcal A$-algebra isomorphism
$$\mathcal A\otimes_{\Z}\Lambda\xrightarrow{\sim} {_\mathcal A} U^-,$$
which sends the Schur functions $s_{\lambda}$ to $IC(\mathcal O_{\lambda})$ for all $\lambda$. By the classical representation theory of the symmetric group, there is an $\mathcal A$-algebra isomorphism
$$\mathcal A\otimes_{\Z}(\bigoplus_n K_0(S_n)) \xrightarrow{\sim} \mathcal A\otimes_{\Z}\Lambda,$$
which sends the Specht module $S^\lambda$ to $s_{\lambda}$. Furthermore, we have an $\mathcal A$-algebra involution $\omega$ of $\mathcal A\otimes_{\Z}(\bigoplus_n K_0(S_n))$ sending $S^\lambda$ to $S^{\widetilde{\lambda}}$.
\end{proof}
\end{lemma}
\vskip 3mm

\begin{example}The above lemma has the following straightforward application. By \cite{C08}, we have for each $n$, the Springer functor  $$\Spr_n\otimes_{S_n}-: \K[S_n]\text{-}\text{mod}\rightarrow P_{GL_{ni}}(E_{ni}^{nil}), $$
which is an equivalence of categories, mapping  irreducible modules to   $IC(\mathcal O_{\lambda}), \lambda\vdash n$.  In particular, we have $$\Spr_n\otimes_{S_n}(S^{(n)})=IC(\mathcal O_{(n)}),\quad \Spr_n\otimes_{S_n}(\K[S_n])=\Spr_n.$$ According to \cite[Theorem 1.3]{C08}, these functors induce  an algebra isomorphism
$$\bigoplus_{n}(\Spr_n\otimes_{S_n}-): \mathcal A\otimes_{\Z}(\bigoplus_n K_0(S_n)) \xrightarrow{\sim} {_\mathcal A} U^-.$$
Thus, by Lemma \ref{springer}, we obtain $\Spr_n\otimes_{S_n}(S^\lambda)=IC(\mathcal O_{\lambda})$ for all $\lambda\vdash n$.
\end{example}
\vskip 3mm

Let $I= I^0=\{i\}$. Then, $K_0(R)=\bigoplus_{n}K_0(R(ni))$. There is an obvious $\mathcal A$-algebra isomorphism
$$\Theta:\mathcal A\otimes_{\Z}(\bigoplus_n K_0(S_n)) \xrightarrow{\sim}K_0(R),$$
which sends the irreducible modules to their gr-projective covers. We thus obtain an  $\mathcal A$-algebra isomorphism $\Theta'=\Gamma^{-1}\Theta:\mathcal A\otimes_{\Z}(\bigoplus_n K_0(S_n))\rightarrow {_\mathcal A} U^-$  with $\Theta'(S^{(n)})=F_{in}$. By Lemma \ref{springer}, we conclude that  $\Theta'(S^\lambda)=IC(\mathcal O_{\widetilde{\lambda}})$ for all $\lambda\vdash n$. Therefore, $\Gamma^{-1}:K_0(R)\rightarrow {_\mathcal A} U^-$ sends the self-dual indecomposable projective modules   to the canonical bases.

\vskip 6mm
\section{\textbf{Categorification of irreducible highest weight module in Jordan quiver case}}
\vskip 3mm
We show in  Jordan quiver case that the cyclotomic KLR-algebras  provide a categorification of the irreducible highest weight $U_q(\g)$-modules.
\vskip 3mm
\subsection{The algebra $U_q(\g)$ and its irreducible highest weight modules}\
\vskip 3mm
Given a Borcherds-Cartan datum $(I, A,\cdot)$.  We set
\begin{itemize}
\item[(a)] $P^\vee=(\bigoplus_{i\in I}\Z h_i)\oplus (\bigoplus_{i\in I}\Z d_i)$, a free abelian group,  the {\it dual weight lattice},
\item[(b)] $\h=\Q\otimes_{\Z}P^\vee$, the {\it Cartan subalgebra},
\item[(c)] $P=\{\lambda \in \h^*\mid \lambda(P^\vee)\subseteq \Z\}$, the {\it weight lattice},
\item[(d)] $\Pi^\vee=\{h_i\in P^\vee\mid i\in I\}$, the set of {\it simple coroots},
\item[(e)] $\Pi=\{\alpha_i\in P \mid i\in I\}$, the set of {\it simple roots}, which is linearly independent over $\Q$ and satisfies
$$\alpha_j(h_i)=a_{ij},\ \alpha_j(d_i)=\delta_{ij}\ \ \text{for all}\ i,j \in I,$$
\item[(f)] for each $i\in I$, there is a $\Lambda_i\in P$,  called the {\it fundamental weight}, defined by
$$\Lambda_i(h_j)=\delta_{ij},\ \Lambda_i(d_j)=0 \ \ \text{for all} \ i,j\in I.$$
\end{itemize}
\vskip 3mm

Let $P^+=\{\Lambda \in P \mid \Lambda(h_i)\geq 0 \  \text{for all} \ i \in I \}$ be the set of dominant integral weights. The free abelian group $Q=\bigoplus_{i \in I} {\Z \alpha_i}$ is called the root lattice. We identify $\Pi$ with $I$ and identify the positive root lattice $Q_+=\bigoplus_{i \in I}{\N\alpha_i}$ with $\N[I]$.

\vskip 3mm

We extend the bilinear form `$\cdot$' to  a non-degenerate symmetric bilinear form $( \ , \ )$ on $\h^*$ satisfying
$$(\alpha_i, \lambda)=r_i\lambda(h_i), \ (\Lambda_i,\lambda)=r_i\lambda(d_i) \ \ \text{for any} \ \lambda \in \h^* \ \text{and}\ i\in I,$$
and therefore we have
$(\alpha_i,\alpha_j)=i\cdot j= r_ia_{ij}=r_ja_{ji}\ \ \text{for all} \ i,j \in I$.




\vskip 3mm

For this extended  datum. We denote by $\widehat {U}$  the $\Q(q)$-algebra generated by the elements $q^h$ $(h\in P^{\vee})$ and $E_{i\ell},
F_{i\ell}$ $((i,\ell) \in I^{\infty})$, satisfying
\begin{equation*}
\begin{aligned}
& q^0= 1,\quad q^hq^{h'}=q^{h+h'} \ \ \text{for} \ h,h' \in P^{\vee} \\
& q^h E_{j\ell}q^{-h} = q^{\ell\alpha_j(h)} E_{j\ell}, \ \ q^h F_{j\ell}q^{-h} = q^{-\ell\alpha_j(h)} F_{j\ell}\ \ \text{for} \ h \in P^{\vee}, (j,\ell)\in I^{\infty}, \\
& \sum_{r+s=1-\ell a_{ij}}(-1)^r
{E_i}^{(r)}E_{j\ell}E_i^{(s)}=0 \ \ \text{for} \ i\in
I^{+},\ (j,\ell)\in I^{\infty} \ \text {and} \ i \neq (j,\ell), \\
& \sum_{r+s=1-\ell a_{ij}}(-1)^r
{F_i}^{(r)}F_{j\ell}F_i^{(s)}=0 \ \ \text{for} \ i\in
I^{+},(j,\ell)\in I^{\infty} \ \text {and} \ i \neq (j,\ell), \\
& E_{ik}E_{j\ell}-E_{j\ell}E_{ik} = F_{ik}F_{j\ell}-F_{j\ell}F_{ik} =0 \ \ \text{for} \ a_{ij}=0,
\end{aligned}
\end{equation*}
which is $Q$-graded by assigning $|q^h|=0$, $|E_{i\ell}|= \ell \alpha_{i}$ and $|F_{i\ell}|= - \ell \alpha_{i}$.

\vskip 3mm

The algebra $\widehat{U}$ is endowed with a co-multiplication
$\Delta\colon \widehat{U} \rightarrow \widehat{U} \otimes \widehat{U}$
given by
$$
\begin{aligned}
& \Delta(q^h) = q^h \otimes q^h, \\
& \Delta(E_{i\ell}) = \sum_{m+n=\ell} q_{(i)}^{mn}E_{im}K_{i}^{n}\otimes E_{in}, \\
& \Delta(F_{i\ell}) = \sum_{m+n=\ell} q_{(i)}^{-mn}F_{im}\otimes K_{i}^{-m}F_{in},
\end{aligned}
$$
where $K_i=q_i^{h_i}$ $(i \in I)$.

\vskip 3mm
Let $\omega\colon \widehat{U}\rightarrow\widehat{U}$  be the $\Q(q)$-algebra involution given by
$$\omega(q^h)=q^{-h},\ \omega(E_{i\ell})=F_{i\ell},\ \omega(F_{i\ell})=E_{i\ell}\ \ \text{for}\ h \in P^{\vee},\ (i,\ell)\in I^{\infty}.$$
Let $\widehat{U}^{+}$ (resp. $\widehat{U}^-$, resp.   $\widehat{U}^{\geq 0}$) be the subalgebra of $\widehat{U}$ generated by $E_{i\ell}$ $((i,\ell) \in I^{\infty})$ (resp. $F_{i\ell}$ $((i,\ell) \in I^{\infty})$, resp. $E_{i\ell}$ $((i,\ell) \in I^{\infty})$ and $q^h$ $(h\in P^{\vee})$). We  identify $\widehat{U}^{-}$ with the $U^-$  in Definition \ref{U-}, and define a symmetric bilinear form $\{ \ , \ \}$ on $\widehat{U}^{\geq 0}$ by setting
$$
\begin{aligned}
&  \{x,y\}=\{\omega(x),\omega(y)\}\ \ \text{for}\ x,y\in \widehat{U}^{+},\\
& \{q^h,1\}=1,\ \{q^h,E_{i\ell}\}=0, \ \ \{q^h,K_j\}=q^{\alpha_j(h)}.
\end{aligned}
$$

\vskip 3mm
\begin{definition} By the Drinfeld double process, we define the  quantum Borcherds-Bozec algebra  $U =U_q(\g)$   as the quotient of $\widehat{U}$ by the relations
\begin{equation}\label{drinfeld}
\sum\{a_{(1)},b_{(2)}\}\omega(b_{(1)})a_{(2)}=\sum\{a_{(2)},b_{(1)}\}a_{(1)}\omega(b_{(2)})\ \ \text{for all}\ a,b \in \widehat{U}^{\geq0}.
\end{equation}
Here we use the Sweedler's notation and write
$\Delta(x)=\sum x_{(1)}\otimes x_{(2)}$.
\end{definition}

 The subalgebra $U^-$  of $U$ generated by   $F_{i\ell}, (i,\ell)\in I^{\infty}$  coincides with  Definition \ref{U-}.

\vskip 3mm


 Let $\Lambda\in P^+$. The irreducible highest weight $U_q(\g)$-module $V(\Lambda)$ is given by
\begin{equation}\label{irr}
\begin{aligned}
 &V(\Lambda)\cong U\bigg/\left (\sum_{(i,\ell)\in I^{\infty}}UE_{i\ell}+\sum_{h\in P^\vee}U(q^h-q^{\Lambda(h)})+\sum_{i\in I^{+}}U F_i^{\Lambda(h_i)+1}+
\sum_{i\in I^{\leq 0} \ \text{with}\ \Lambda(h_i)=0;\ \ell\geq 1}UF_{i\ell}\right) \\
 &\phantom{V(\Lambda)}\cong U^- \bigg/ 
\left(\sum_{i\in I^{+}}U^- F_i^{\Lambda(h_i)+1}+
\sum_{i\in I^{\leq 0} \ \text{with}\ \Lambda(h_i)=0; \ \ell\geq 1}U^-F_{i\ell}\right).
\end{aligned}
\end{equation}
\vskip 5mm

\subsection{Jordan quiver case}\
\vskip 3mm
Throughout this section, we assume that $I=I^0=\{i\}$ and $U$ is the quantum Borcherds-Bozec associated to $I$.  By \cite[Appendix]{FKKT2021}, if we define $\{\alpha_p\}_{p\geq 1}$ inductively:
$$
\begin{aligned} \alpha_p=\nu_p(K_i^{-p}-K_i^{p})-\nu_1K_i\alpha_{p-1}-\nu_2K_i^{2}\alpha_{p-2}-\cdots-\nu_{p-1}K_i^{p-1}\alpha_1,
\end{aligned}
$$
where $\alpha_1=\nu_1(K_i-K_i^{-1})$, then for any $\ell,t\geq 1$, the equation (\ref{drinfeld}) yields
\begin{equation}\label{comm}[E_{i\ell},F_{it}]=\sum_{p=0}^{\min{\{\ell,t\}}}\alpha_p F_{i,t-p}E_{i,\ell-p}.\end{equation}

\vskip 3mm

Define the functors
$$\begin{aligned}
& \F_{i\ell}\colon R(ni)\text{-}\Mod\rightarrow R((n+\ell)i)\text{-}\Mod,\  M\mapsto (R((n+\ell)i) 1_{ni}\otimes e_{i,\ell})\otimes_{R(ni)}M,\\
& \E_{i\ell}\colon R(ni)\text{-}\Mod\rightarrow R((n-\ell)i)\text{-}\Mod,\ M\mapsto 1_{(n-\ell)i}\otimes e_{i,\ell}M.
\end{aligned}$$

\vskip 3mm
\begin{lemma}\label{LLL8}\
{\it Let $\ell,t\geq 1$. We have the following natural isomorphisms
$$
\E_{i\ell}\F_{it}\simeq \bigoplus_{p=0}^{\min{\{\ell,t\}}}\F_{i,t-p}\E_{i,\ell-p}\otimes  Z_{p},
$$
where $Z_{p}$ is the algebra of symmetric polynomials in $p$  indeterminates, each of degree $2r_i$.
}
\begin{proof}
We prove the case where $t=\ell$  only. The other cases are similar. For simplicity, we omit the symbol  ``$i$".
Assume first that $n\geq \ell$ and denote by $D_{n,\ell}$ (resp. $D^{-1}_{n,\ell}$)  the set of minimal length left (resp. right) $S_{n}\times  S_{\ell}$-coset representatives in $S_{n+\ell}$. Then $D_{n,\ell}\cap D^{-1}_{n,\ell}=\{v_0,\dots,v_\ell\}$ is the set of minimal length $S_{n}\times  S_{\ell}$-double coset representatives, where $v_k$ can be expressed graphically

$$\coset$$

 Note that

\begin{equation}\label{4.1}1_n\otimes e_\ell\cdot R(n+\ell)\cdot 1_n\otimes e_\ell=\sum_{u\in D_{n,\ell}} 1_n\otimes e_\ell\cdot u \cdot R(n)\otimes R(\ell)\cdot 1_n\otimes e_\ell.\end{equation}
Any $u\in D_{n,\ell}$ can be decomposed into $u=v\otimes v'\cdot v_k$ for some $0\leq k\leq \ell$, $v\in D_{n-k,k}$ and $v'\in D_{k,\ell-k}$. Since $e_\ell\cdot v'=e_\ell$, we have
\begin{equation}\label{4.2}\begin{aligned} & 1_n\otimes e_\ell\cdot R(n+\ell)\cdot 1_n\otimes e_\ell=\sum_{k=0}^{\ell}\sum_{v\in D_{n-k,k}} 1_n\otimes e_\ell\cdot v\cdot v_k \cdot R(n)\otimes R(\ell)\cdot 1_n\otimes e_\ell
\\&\phantom{1_n\otimes e_\ell \cdot R(n+\ell)\cdot 1_n\otimes e_\ell}=\bigoplus_{k=0}^{\ell}\sum_{v\in D_{n-k,k}} v\otimes e_\ell\cdot v_k \cdot R(n)\otimes R(k)\otimes R(\ell-k)\cdot 1_n\otimes e_\ell
\end{aligned}\end{equation}
On the other hand, $$(R(n)\cdot 1_{n-k}\otimes e_k)\otimes_{R(n-k)}(1_{n-k}\otimes e_k\cdot R(n))=\sum_{v\in D_{n-k,k}}v\cdot (1_{n-k}\otimes R(k)e_k)\otimes_{R(n-k)}(1_{n-k}\otimes e_k\cdot R(n)).$$
Since $(1_{n-k}\otimes R(k)e_k)\cdot v_k=v_k\cdot (1_{n}\otimes R(k)e_k)$ and $v_k\cdot 1_{n-k}\otimes e_k\cdot R(n)=1_n\otimes e_{k}\cdot v_k \cdot R(n)$, we see that
$$\begin{aligned} & 1_n\otimes e_\ell\cdot \Big(v\cdot (1_{n-k}\otimes R(k)e_k)\Big)\cdot v_k\cdot 1_{n+k}\otimes e_{\ell-k}R(\ell-k)e_{\ell-k}\cdot\Big((1_{n-k}\otimes e_k\cdot R(n))\Big)\cdot 1_n\otimes e_\ell \\ &= 1_n\otimes e_\ell \cdot \Big(v\cdot v_k\cdot R(n)\otimes R(k)e_k\otimes e_{\ell-k}R(\ell-k)e_{\ell-k}\Big)\cdot 1_n\otimes e_\ell. \\
& = v\otimes e_\ell\cdot v_k \cdot R(n)\otimes R(k)\otimes R(\ell-k)\cdot 1_n\otimes e_\ell{\color{red}.}\end{aligned} $$
Hence for each $z\in e_{\ell-k}R(\ell-k)e_{\ell-k}$, the map
$$(R(n)\cdot 1_{n-k}\otimes e_k)\otimes_{R(n-k)}(1_{n-k}\otimes e_k\cdot R(n))\rightarrow 1_n\otimes e_\ell \cdot R(n+\ell)\cdot 1_n\otimes e_\ell$$
$$x\otimes y\mapsto 1_n\otimes e_\ell\cdot x\cdot 1_{n+k}\otimes z\cdot y\cdot 1_n\otimes e_\ell$$
is an injective $(R(n),R(n))$-bimodule homomorphism. Since $e_{\ell-k}R(\ell-k)e_{\ell-k}\cong Z_{\ell-k}$, we have proved that
\begin{equation}\label{4.3}\E_{\ell}\F_{\ell}\simeq \bigoplus_{k=0}^{\ell}\F_{k}\E_{k}\otimes  Z_{\ell-k}\end{equation} on $R(n)$. If $n<\ell$, then the direct sum in (\ref{4.2}) ranges  from $k=0$ to $n$, while the right hand side of (\ref{4.3}) only makes sense for $k\leq n$.
\end{proof}
\end{lemma}
\vskip 3mm

Choose   $\Lambda\in P^+$ and set $a=\Lambda(h_i)\geq 0$. We define the cyclotomic algebra $R^{\Lambda}(n)$ to be the quotient of $R(n)$ by the two sided ideal generated by $x_{1}^{a}$, and form
$$R^\Lambda=\bigoplus_{n\geq 0}R^{\Lambda}(n), \ \ \ K_0(R^\Lambda)=\bigoplus_{n\geq 0}K_0(R^\Lambda(n)).$$
If $a=0$, then $R^\Lambda=R^\Lambda(0)=\K$ and $V(\Lambda)$ is the one dimensional trivial module by (\ref{irr}). So we  assume that  $a>0$ in the following.
\vskip 3mm
Note that $R^\Lambda(n)$ has a basis $\{x_1^{r_1}\cdots x_n^{r_n}\tau_\omega\mid \omega\in S_n,0\leq r_1,\dots,r_n<a\}$. Define  the functors
$$\begin{aligned}
& \F^\Lambda_{i\ell}\colon R^\Lambda(n)\text{-}\Mod\rightarrow R^\Lambda(n+\ell)\text{-}\Mod,\  \ M\mapsto (R^\Lambda(n+\ell) 1_{n}\otimes e_{\ell})\otimes_{R^\Lambda(n)}M,\\
& \E^\Lambda_{i\ell}\colon R^\Lambda(n)\text{-}\Mod\rightarrow R^\Lambda(n-\ell)\text{-}\Mod,\ \ M\mapsto 1_{n-\ell}\otimes e_{\ell}M.
\end{aligned}$$

\vskip 3mm
Similar to Lemma \ref{LLL8}, for  $\ell,t\geq 1$, we have the following natural isomorphisms
\begin{equation}\label{4.4}
\E^\Lambda_{i\ell}\F^\Lambda_{it}\simeq \bigoplus_{p=0}^{\min{\{\ell,t\}}}\F^\Lambda_{i,t-p}\E^\Lambda_{i,\ell-p}\otimes  Z^\Lambda_{p},
\end{equation}
where $Z_{p}^\Lambda={\left(\K[x_1,\dots,x_p]/(x_1^a,\dots,x_p^a)\right)}^{S_p}$, i.e., the symmetric polynomials in $x_1,\dots,x_p$ such that   no $x_k^m$ $(m\geq a)$ appears. Thus $Z_p^\Lambda$ is determined by  all   partitions $\lambda$ with $\ell(\lambda)\leq p$ and $\lambda_1\leq a-1$. We know that the generating function for such partitions is

$${\begin{bmatrix} a+p-1 \\ p \end{bmatrix}}= \frac{(1-q^a)(1-q^{a+1})\cdots(1-q^{a+p-1})}{(1-q)(1-q^2)\cdots(1-q^p)},$$
and therefore
$$\gdim \, Z_{p}^\Lambda=\frac{(1-q_i^{2a})(1-q_i^{2(a+1)})\cdots(1-q_i^{2(a+p-1)})}{(1-q_i^2)(1-q_i^4)\cdots(1-q_i^{2p})}:=\beta_p.$$
\vskip 3mm
\begin{lemma}
{\it Let $\nu_k=1/(1-q_i^2)(1-q_i^4)\cdots(1-q_i^{2k})$ for any $k\geq 1$. Then for any $p\geq 1$, \begin{equation}\label{Gauss1}\beta_p=\nu_p(1-q_i^{2pa})-\nu_1q_i^{2a}\beta_{p-1}-\nu_2q_i^{4a}\beta_{p-2}-\cdots-\nu_{p-1}q_i^{2(p-1)a}\beta_{1}.\end{equation}
}
\begin{proof}
Using the notations
$$\begin{bmatrix} n \\ m \end{bmatrix} =\frac{(1-q^n)(1-q^{n-1})\cdots(1-q^{n-m+1})}{(1-q)(1-q^2)\cdots(1-q^m)}$$
for $n\geq m\geq 1$ and $(x;q)_n=(1-x)(1-xq)\cdots (1-xq^{n-1})$ for $n\geq 1$. We have
\begin{equation}\label{Gauss}\begin{bmatrix} n+1 \\ m \end{bmatrix}= q^m\begin{bmatrix} n \\ m \end{bmatrix}+\begin{bmatrix} n \\ m-1 \end{bmatrix},\quad (xq^m;q)_{n-m}=\frac{(x;q)_n}{(x;q)_m}.\end{equation}

To show the identity in the lemma,  it is enough to show the following
$$(1-q^{pa})=\sum_{k=0}^{p-1}q^{ka}\begin{bmatrix} p \\ k \end{bmatrix}(q^a;q)_{p-k},$$
which can be proved easily by an induction on $p$ and using  (\ref{Gauss}).
\end{proof}
\end{lemma}

\vskip 3mm
Let $\alpha_k=q_i^{-ka}\beta_k$ for all $k\geq 1$. By (\ref{Gauss1}), for $p\geq 1$, we have
$$q_i^{-pa}\beta_p=\nu_p(q_i^{-pa}-q_i^{pa})-\sum_{k=1}^{p-1}q_i^{-(p-k)a}\nu_kq_i^{ka}\beta_{p-k}.$$
Thus,
$$\alpha_p=\nu_p(q_i^{-pa}-q_i^{pa})-\sum_{k=1}^{p-1}\nu_kq_i^{ka}\alpha_{p-k}.$$
\vskip 3mm
Define the functors $ E_{i\ell}^\Lambda,  F_{i\ell}^\Lambda, K_i$ on  $K_0(R^\Lambda)$  by
$$E_{i\ell}^\Lambda=\mathcal E_{i\ell}^\Lambda,\ \ \  F_{i\ell}^\Lambda= q_i^{-\ell a}\mathcal F_{i\ell}^\Lambda, \ \ \ K_i=q_i^a.$$
Then (\ref{4.4}) gives
$$E^\Lambda_{i\ell}F^\Lambda_{it}= \sum_{p=0}^{\min{\{\ell,t\}}}q_i^{-pa}\beta_pF^\Lambda_{i,t-p}E^\Lambda_{i,\ell-p}= \sum_{p=0}^{\min{\{\ell,t\}}}\alpha_pF^\Lambda_{i,t-p}E^\Lambda_{i,\ell-p},$$
where $\alpha_p=\nu_p(K_i^{-p}-K_i^{p})-\sum_{k=1}^{p-1}\nu_kK_i^{k}\alpha_{p-k}$.
\vskip 3mm
Let $K_0(R^\Lambda)_{\Q(q)}=\Q(q)\otimes_{\Z[q,q^{-1}]}K_0(R^\Lambda)$. Then by (\ref{comm}), the $K_0(R^\Lambda)_{\Q(q)}$  is a weight $U$-modules (the weight spaces are $K_0(R^\Lambda)_{\Lambda-ni}=K_0(R^\Lambda(n))$) with the action of $E_{i\ell}$ (resp. $F_{i\ell}$) by $ E_{i\ell}^\Lambda$ (resp. $F_{i\ell}^\Lambda$). The $\Z[q,q^{-1}]$-linear map
$$\varphi\colon K_0(R)\rightarrow K_0(R^\Lambda), \ \ [P]\mapsto R^\Lambda(n) \otimes_{R(n)}[P]$$
is an isomorphism. For $[P]\in K_0(R(n))$, we have
$$\begin{aligned}&\varphi(F_{i\ell}[P])=\varphi(\Ind_{n,\ell}^{n+\ell} P\otimes R(\ell)e_\ell)\\
& \phantom{\varphi(F_{i\ell}[P]}=R^\Lambda(n+\ell) \otimes_{R(n+\ell)}R(n+\ell)\otimes_{R(n)\otimes R(\ell)}P\otimes  R(\ell)e_\ell\\
& \phantom{\varphi(F_{i\ell}[P]}=R^\Lambda(n+\ell)\otimes_{R(n)\otimes R(\ell)}P\otimes  R(\ell)e_\ell\\
& \phantom{\varphi(F_{i\ell}[P]}=(R^\Lambda(n+\ell)1_n\otimes e_\ell)\otimes_{R(n)}P\\
& \phantom{\varphi(F_{i\ell}[P]}=\F_{i\ell}^\Lambda \varphi([P]).
\end{aligned}$$
It follows that $\varphi$ is $U^-$-linear and  $K_0(R^\Lambda)$ is generated by $1_\Lambda$,  the trivial module over $R^\Lambda(0)$.
Hence $K_0(R^\Lambda)_{\Q(q)}$ is isomorphic to the irreducible highest weight module $V(\Lambda)$ given in (\ref{irr}), which can be identified with $U^-$ as $U^-$-modules.

\vskip 3mm
\begin{theorem}{\it
If $I=I^0=\{i\}$, then $K_0(R^\Lambda)_{\Q(q)}$ is isomorphic to the irreducible highest weight module $V(\Lambda)$ for each $\Lambda\in P^+$.
}\end{theorem}

\vskip 6mm

\appendix
\section{The KLR-algebras of  $\mathcal K^1$}
\vskip 3mm
Recall that $\mathcal K^1$ is the subalgebra of $U^-$ generated by  $F_i$ for $i\in I^+\cup I^-$,  and $F_{i\ell}$ for $i\in I^{ 0}, \ell>0$. Fix $\alpha=\sum_{i\in I}\alpha_ii\in \N[I]$ with $\text{ht}{(\alpha)}:=\sum_{i\in I}\alpha_i=n$. Let $\text{Seq}(\alpha)$ be the set of all sequences $\ii=i_1i_2\dots i_n$ in $I$ such that $\alpha=i_1+i_2\cdots +i_n$.

\vskip 3mm

We define the Khovanov-Lauda-Rouquier algebra
$\mathcal R(\alpha)$ associated to a given Borcherds-Cartan datum $(I,A,\cdot)$  to be the $\K$-algebra with the homogeneous generators:

$$1_{\ii}=\genO{i_1}{i_k}{i_n} \quad \text{for} \ \ii=i_1i_2\dots i_n\in \text{Seq}(\alpha) \ \text{with} \ \text{deg}(1_{\ii})=0,$$
\vskip 2mm
$$\ \ x_{k,\ii}=\genX{i_1}{i_k}{i_n} \quad \text{for} \ \ii\in \text{Seq}(\alpha),1\leq k\leq n \ \text{with} \ \text{deg}(x_{k,\ii})=2r_{i_k},$$
\vskip 2mm
$$\tau_{k,\ii}=\genTT{i_1}{i_k}{i_{k+1}}{i_n} \quad \text{for} \ 1\leq k\leq n-1 \ \text{with} \ \text{deg}(\tau_{k,\ii})=-i_k\cdot i_{k+1}.$$

\vskip 3mm
\noindent subject to the following local relations:

\begin{align}
 \dcrossS{i}{j} \  =  \  \begin{cases} \quad \quad \quad \quad \quad \  0 & \text{ if } i= j\in I^+, \\  \\ \   \Big( \dcrossLS{-\frac{a_{ii}}{2}}{i}{i} \ + \ \dcrossRS{-\frac{a_{ii}}{2}}{i}{i}\Big)^2  & \text{ if } i= j \ \text{and} \  i\cdot i < 0,
                    \\ \\
                    \quad  \quad \quad \quad \ \ \dcrossAS{i}{j} & \text{ if } \ i\cdot j=0, \\
                    \\
                    \   \dcrossLS{-a_{ij}}{i}{j} \ + \ \dcrossRS{-a_{ji}}{i}{j}  & \text{ if } i\ne j \ \text{and} \  i\cdot j < 0,
                  \end{cases}
\end{align}

\begin{equation}
\begin{aligned}
 \LUS{{}}{i}{i} \   -  \ \RDS{{}}{i}{i}\ =\ \dcrossAS{i}{i} \quad\quad\quad \LDS{{}}{i}{i}   \ - \ \RUS{{}}{i}{i}\ =\ \dcrossAS{i}{i}  \quad \text{ if } i \in I^+,
\end{aligned}
\end{equation}

\begin{equation}
\begin{aligned}
 \LUS{{}}{i}{j} \   =  \ \RDS{{}}{i}{j}  \quad\quad\quad \LDS{{}}{i}{j}   \ =  \ \RUS{{}}{i}{j}  \quad \text{ otherwise},
\end{aligned}
\end{equation}

\begin{equation}
 \BraidLS{i}{j}{i} \  -  \ \BraidRS{i}{j}{i} \ =  {\sum_{c=0}^{-a_{ij}-1}} \ \threeDotStrandsS{i}{j}{i} \ \ \text{ if } i\in I^+, i\ne j \ \text{and} \ i\cdot j < 0,
 \end{equation}
 \begin{equation}
\BraidLS{i}{j}{k} \  = \ \BraidRS{i}{j}{k} \quad\quad \text{otherwise}.
\end{equation}

\vskip 3mm
Denote by $K_0(\R)$ ($\R=\bigoplus_{\alpha}\R(\alpha)$) the Grothendieck group of the category of finite generated gr-projective $\R$-modules. As in Section 2, we  endow  $K_0(\R)$ with a twisted bialgebras structure, and subsequently, we obtain a twisted bialgebra isomorphism  $ \mathcal K^1 \xrightarrow{\sim} K_0(\R)$ given by
$$\begin{aligned}& F_i^{(n)}\mapsto [P_{i^{(n)}}]\quad \text{for}\ i\in I^+,n\geq 0,\\
 & F_{i\ell}\mapsto [P_{i,\ell}]\qquad  \text{for}\ i\in I^0, \ell\geq 1,\\
 & F_i\mapsto [P_{i}]\quad\quad \text{for}\ i\in I^-. \end{aligned}$$
We also conjecture that the indecomposable projective modules of $K_0(\R)$ are mapped to the canonical basis   $\bigsqcup_{\alpha\in\N[I]}\mathcal P_{\alpha}^1$ of $\mathcal K^1$.

\vskip 10mm

\bibliographystyle{amsplain}

\end{document}